\documentclass[a4paper,fleqn]{cas-sc}

\usepackage{float}
\usepackage{algorithm, algpseudocode}
\usepackage{amsthm}
\newtheorem{proposition}{Proposition}
\usepackage{mathtools}
\DeclareMathOperator*{\argmin}{argmin}

\usepackage[numbers,sort&compress]{natbib}

\def\tsc#1{\csdef{#1}{\textsc{\lowercase{#1}}\xspace}}
\tsc{WGM}
\tsc{QE}

\begin{document}
\let\WriteBookmarks\relax
\def\floatpagepagefraction{1}
\def\textpagefraction{.001}

\shorttitle{}    

\shortauthors{}  

\title [mode = title]{Distributed Stochastic Model Predictive Control with Temporal Aggregation for the Joint Dispatch of Cascaded Hydropower and Renewables}  



%

\author[1,2]{Luca Santosuosso}[orcid=0000-0002-7865-3883]
\cormark[1]
\ead{luca.santosuosso@tugraz.at}

\author[1,2]{Sonja Wogrin}[orcid=0000-0002-3889-7197]
\ead{wogrin@tugraz.at}

\affiliation[1]{
            organization={Institute of Electricity Economics and Energy Innovation, Graz University of Technology},
            addressline={Inffeldgasse 18},
            city={Graz},
            postcode={8010},
            country={Austria}
}
\affiliation[2]{
    organization={Research Center ENERGETIC},
    addressline={Rechbauerstraße 12},
    city={Graz},
    postcode={8010},
    country={Austria}
}

\cortext[1]{Corresponding author}



\begin{abstract}
This paper addresses the real-time energy dispatch of a hybrid system comprising cascaded run-of-the-river hydropower plants, wind, and solar photovoltaic units, operated under uncertainty in water inflows and renewable power generation.
Traditional scenario-based stochastic model predictive control (MPC) schemes suffer from severe computational limitations due to the high dimensionality induced by both the temporal and scenario dimensions of the dispatch problem, 
as well as the inherent nonconvexities associated with cascaded hydropower dynamics.
To overcome these challenges, we propose a novel control scheme that seamlessly integrates time series aggregation (TSA), distributed optimization, and stochastic MPC.
The resulting temporally aggregated distributed stochastic MPC scheme simultaneously reduces the temporal dimension of the dispatch problem via TSA and decomposes it across scenarios through distributed optimization.
Our main theoretical result establishes a formal performance guarantee for the proposed controller,
enabling a rigorous quantification of its solution accuracy at every MPC iteration.
Numerical results based on a real-world case study show the effectiveness of the proposed controller,
achieving up to 74\% reduction in computational effort relative to the full-scale centralized counterpart when the required solution accuracy is at least 99\%,
and up to 85\% when the accuracy requirement is relaxed to 95\%.
Notably, the proposed controller not only significantly enhances computational efficiency relative to the traditional full-scale centralized counterpart, but more importantly restores computational tractability, whereas the traditional controller fails to solve the dispatch problem within the prescribed time limit for computing control actions.
\end{abstract}


\maketitle

\section{Introduction}
\label{sec:Introduction}

The inherent stochasticity of variable renewable energy sources (vRES), such as wind and solar photovoltaic, has spurred increasing interest in their joint dispatch with controllable units \citep{jing2024benefit}.
Cascaded hydropower plants, which harness the water potential at multiple points along a river, have proven particularly effective for this purpose \citep{apostolopoulou2018robust}.
By combining clean power generation, fast ramping capabilities, and the storage capacity of water reservoirs, such systems, when coupled with vRES, form hybrid configurations capable of mitigating unforeseen power fluctuations internally,
while actively participating in energy trading \citep{ma2022decentralized}, peak shaving \citep{wu2024short}, and the provision of ancillary services \citep{santosuosso2023economic}.

The efficient operation of cascaded hydropower-vRES (CH-vRES) hybrid systems typically involves multiple decision-making stages, ranging from long-term planning to short-term scheduling and fast, reactive real-time control \citep{su2025real}.
This paper focuses on the latter stage of this sequential decision-making process.

The real-time control of CH-vRES systems typically operates at sub-hourly resolutions \citep{hamann2014real}.
Among the available control strategies, model predictive control (MPC) is widely adopted due to its ability to explicitly incorporate operational constraints and control objectives,
while enabling iterative decision-making through the continuous updating of control actions as the uncertainties are progressively revealed \citep{ye2023real}.

When MPC is applied to CH-vRES systems, accurately representing the underlying nonconvex physical dynamics typically leads to a mixed-integer nonlinear programming (MINLP) formulation \citep{karamyar2025scalable}, for which efficient, general-purpose off-the-shelf solvers remain unavailable.
Consequently, a substantial body of research has focused on developing surrogate formulations, spanning both linear \citep{liu2026day} and nonlinear \citep{lu2024medium} optimization models, to strike a pragmatic balance between modeling fidelity and computational tractability.
Among these formulations, mixed-integer programming (MIP) models \citep{zhao2025hydro}, particularly in the form of mixed-integer linear programming (MILP) models \citep{zhang2025comparative} and mixed-integer quadratic programming (MIQP) models \citep{liao2021daily}, have emerged as particularly effective.
Nevertheless, while such formulations may be tractable in static or long-term optimization settings,
their implementation at high temporal resolutions within real-time dispatch frameworks often remains computationally prohibitive \citep{wang2026deep}.

The difficulties associated with handling the nonconvex physical dynamics inherent in CH-vRES operation are further compounded by multiple sources of uncertainty, most notably hydrological inflows and variable renewable energy generation \citep{qiu2020stochastic}.
Nevertheless, the integration of uncertainty-aware decision-making within MPC schemes remains essential for the practical operation of CH-vRES systems \citep{zhang2022long}.
In this context, robust optimization \citep{zhou2021100}, chance-constrained programming \citep{zhang2022chance}, and stochastic programming \citep{shi2024stochastic} constitute the predominant methodological paradigms.

Robust MPC assumes that uncertainties are confined within predefined uncertainty sets; however, this assumption may be overly restrictive or even impractical for inherently unbounded or poorly characterized uncertainty sources, often leading to overly conservative control policies \citep{saltik2018outlook}.
Chance-constrained MPC, while offering a probabilistic treatment of constraint satisfaction, typically requires intricate reformulations that increase both the analytical and computational complexity of the resulting control scheme \citep{lyons2012chance}.

In contrast, scenario-based stochastic MPC is often preferred in practical implementations, as it provides a flexible and intuitive modeling framework in which the uncertainty domain is discretized into a finite set of possible realizations (or \textit{scenarios}) with associated probabilities \citep{jeong2023implementation}.
However, accurately capturing multiple sources of uncertainty in CH-vRES operation typically requires a large number of scenarios \citep{zhang2019coordinated}.
When combined with the intrinsic nonconvexity of CH-vRES dynamics, this leads to large-scale stochastic nonconvex optimization models that are exceedingly difficult to solve within the stringent time constraints of real-time MPC \citep{velarde2019scenario}.

To tackle this computational challenge, prior research has employed mathematical decomposition methods,
such as Benders decomposition \citep{lopez2018stochastic} and augmented Lagrangian relaxation methods \citep{11270217},
notably the alternating direction method of multipliers (ADMM),
to decompose the centralized MPC scheme into subproblems that are solved iteratively in parallel while seeking convergence towards global optimality, thereby yielding a distributed MPC scheme \citep{santosuosso2024distributed}.

A key limitation of classical decomposition methods is that convergence to a globally optimal solution is generally guaranteed only when the resulting subproblems satisfy convexity assumptions \citep{molzahn2017survey}.
This requirement underpins, for instance, the theoretical convergence guarantees of classical Benders decomposition \citep{rahmaniani2017benders} and ADMM \citep{boyd2010distributed}.
In the context of CH-vRES dispatch problems, however, the presence of nonconvex physical dynamics typically violates these assumptions, thereby limiting the practical applicability of such classical decomposition methods to heuristic solution algorithms in general settings.

Furthermore, existing decomposition-based approaches for reducing the computational complexity of hydropower dispatch under uncertainty predominantly rely on scenario decomposition \citep{moiseeva2017strategic}.
This technique typically involves relaxing the non-anticipativity constraints of the stochastic programming formulation, which in this context is often cast as a two-stage stochastic program \citep{rodriguez2021accelerating}, thereby yielding a decomposed structure in which scenario-wise subproblems (e.g., one per scenario) can be solved in parallel \citep{ge2019multiple}.
While this approach effectively enhances scalability with respect to the number of scenarios in stochastic MPC schemes, it does not explicitly address the temporal dimension of the problem, which represents another key source of computational complexity.
Although decomposition methods could, in principle, also be used to decouple the intertemporal constraints governing CH-vRES operations into subproblems defined over individual time periods, 
such a temporal decomposition would result in a distributed MPC scheme comprising an impractically large number of subproblems, thereby limiting its operational viability.

Alternatively, time series aggregation (TSA) has proven effective for this purpose \citep{hoffmann2020review}.
By condensing the input time series of the MPC scheme (e.g., hydrological inflows and vRES generation scenarios) into a reduced set of representative time periods,
TSA yields a lower-dimensional (or \textit{temporally aggregated}) optimization model that approximates the original full-scale optimization model while significantly reducing computational complexity \citep{tejada2018enhanced}.
In power system applications, traditional TSA methods typically employ clustering techniques such as k-means \citep{liu2017hierarchical} and k-medoids \citep{schutz2018comparison} to identify representative periods,
with the primary objective of preserving the statistical features of the underlying input time series \citep{sarajpoor2023time}.
However, accurately representing the input space of an optimization model does not necessarily guarantee the accuracy of the aggregated model output \citep{wogrin2023time}.
In other words, even if the statistical features of the input time series are accurately captured, this does not generally ensure that the optimal decision variable values or the optimal objective function value of the aggregated model will closely approximate those of the full-scale model \citep{zhang2022model}.
In this sense, traditional TSA-based solution methods often act as purely heuristic methods, lacking a clear measure of the accuracy of the resulting aggregated model \citep{teichgraeber2022time}.

This limitation has motivated growing interest in performance-guaranteed TSA methods \citep{santosuosso2025we},
which aim to derive theoretically validated bounds on the maximum approximation error introduced by the aggregated model relative to its full-scale counterpart \citep{li2022representative}.
When embedded within iterative solution algorithms, these methods enable the construction of aggregated models that provide computable upper and lower bounds on the optimal objective function value of the original full-scale model, thereby allowing an explicit quantification of the resulting optimality gap (i.e., the difference between these bounds) at each iteration \citep{santosuosso2025optimal}.
Despite their significant potential in CH-vRES system applications, where both computational efficiency and solution reliability are critical, no prior study has investigated the use of performance-guaranteed TSA methods in this context.

Notably, establishing performance guarantees for TSA under intertemporal constraints remains particularly challenging, as the aggregated model must preserve consistency with the temporal dynamics of the full-scale model \citep{pineda2018chronological}.
This challenge is further exacerbated in CH-vRES systems, where hydropower reservoirs are subject to complex intertemporal constraints arising from storage dynamics and ramping limits, as well as spatial constraints induced by the cascaded reservoir dynamics.

As the computational complexity of the stochastic CH-vRES dispatch problem grows substantially along both the temporal and scenario dimensions,
the integration of mathematical decomposition and TSA within MPC schemes constitutes a particularly promising, yet largely unexplored, solution paradigm.
To the best of our knowledge, \cite{santosuosso2026stochastic} provides the only prior investigation of such an integrated control scheme.
Specifically, \cite{santosuosso2026stochastic} combines TSA and ADMM within MPC to simultaneously reduce and parallelize the solution of virtual power plant dispatch problems while preserving formal guarantees on the controller's accuracy.
However, the formulation considered in \cite{santosuosso2026stochastic} is convex, and extending this methodology to the highly nonconvex CH-vRES dispatch problems considered here is fundamentally nontrivial,
as the combined presence of nonconvex dynamics and multiple intertemporal couplings directly challenges both components underpinning the performance guarantees established therein.
First, preserving rigorous performance guarantees under TSA requires the temporally aggregated model to retain the relevant theoretical properties despite the nonconvex dynamics and intertemporal constraints inherent in CH-vRES operation.
Second, ADMM is not generally guaranteed to converge for nonconvex optimization problems \cite{boyd2010distributed},
thereby precluding a direct extension of the performance guarantees established for the MPC scheme of \cite{santosuosso2026stochastic}.
Consequently, the methodology developed in \cite{santosuosso2026stochastic} cannot be directly transferred to the problem class addressed in this study.

To date, no existing control scheme has jointly integrated stochastic MPC, mathematical decomposition, and performance-guaranteed TSA
to achieve simultaneous scalability along both the temporal and scenario dimensions of nonconvex CH-vRES dispatch problems
involving multiple intertemporal constraints, such as ramping and reservoir storage constraints,
while enabling a theoretically grounded and transparent characterization of the solution accuracy attained at each MPC iteration.
This paper seeks to bridge this research gap.

The key contributions of this paper are as follows:
\begin{itemize}
\item  We formulate the energy dispatch of CH-vRES hybrid systems as a two-stage stochastic MIQP problem addressed via a conventional full-scale centralized stochastic MPC scheme.
We then propose a novel temporally aggregated centralized stochastic MPC scheme, derived through TSA, to reduce the temporal dimension of the dispatch problem.
We formally demonstrate that, at every MPC iteration, the proposed temporally aggregated centralized stochastic MPC scheme provides a valid lower bound on the optimal objective function value of the original full-scale centralized stochastic MPC scheme.
Importantly, this theoretical result remains valid in the presence of intertemporal, scenario, and spatial coupling constraints, as well as nonconvexities in both stages of the two-stage stochastic dispatch problem.
    
\item To further reduce the computational complexity of the nonconvex CH-vRES dispatch problem along the scenario dimension,
we decompose the temporally aggregated centralized stochastic MPC scheme across scenarios using ADMM, thereby obtaining a temporally aggregated distributed stochastic MPC scheme.
Since ADMM is not generally guaranteed to converge in nonconvex settings, we derive a valid lower bound on the optimal objective function value of the temporally aggregated centralized stochastic MPC scheme.
By construction, this bound also constitutes a valid lower bound on the optimal objective function value of the original full-scale MPC scheme.

\item Building on these theoretical results, we develop an iterative control algorithm that progressively tightens upper and lower bounds on the optimal objective function value of the original full-scale centralized stochastic MPC scheme at every MPC iteration.
The lower bound is obtained by jointly exploiting temporal aggregation and scenario-wise decomposition, whereas the upper bound is constructed through scenario-wise decomposition.
Consequently, the full-scale dispatch problem is never solved at its original resolution.
Crucially, the algorithm generates a feasible solution to the original full-scale dispatch problem at every iteration
and provides a transparent performance guarantee by enabling the continuous evaluation of the achieved optimality gap.

\end{itemize}

Finally, the effectiveness of the proposed controller is evaluated through a case study that mimics the energy dispatch of a real-world CH-vRES hybrid system located along the Rhône river in France \citep{piron2016operating}.

The remainder of the paper is organized as follows.
Section \ref{sec:methodology} presents the proposed methodology,
Section \ref{sec:results} discusses the simulation results,
and Section \ref{sec:conclusions} concludes the study.

\section{Methodology}
\label{sec:methodology}

This section presents the proposed methodology.
Subsection~\ref{subsec:prob_statement} formulates the dispatch problem,
while Subsection~\ref{subsec:SMPC} introduces a conventional full-scale centralized stochastic MPC scheme adopted as the benchmark control method.
Subsection~\ref{subsec:SMPC_TSA} develops a temporally aggregated counterpart of the benchmark controller and presents the main theoretical result of this paper,
namely that the proposed temporally aggregated model always yields a valid lower bound on the optimal objective value of its full-scale counterpart.
Subsequently, Subsection~\ref{subsec:DSMPC_TSA} extends the temporally aggregated controller through scenario decomposition.
Finally, Subsection~\ref{subsec:obj_fun_bounds} presents the proposed control algorithm, which combines stochastic MPC, mathematical decomposition methods, and TSA to simultaneously reduce the temporal and scenario dimensions of the original full-scale centralized stochastic MPC scheme while preserving a formal performance guarantee.

Sets, matrices, and vectors are denoted by boldface symbols.
The cardinality of a set is denoted by $|\cdot|$,
while the Euclidean ($\ell_2$) norm is denoted by $\|\cdot\|_2$.
The zero vector in $\mathbb{R}^m$ is denoted by $\mathbf{0}^m$.
All sets are indexed starting from $0$.

\subsection{Problem Statement}
\label{subsec:prob_statement}

This section introduces the optimal control problem under consideration.
The goal is to determine the optimal dispatch strategy for a hybrid system comprising vRES jointly operated with cascaded hydropower plants.
In particular, we consider a CH-vRES system composed of run-of-the-river hydropower plants integrated with wind and solar power units, consistently with the case study under investigation \citep{santosuosso2025distributed}.
Such hybrid configurations are of growing relevance in several countries where run-of-the-river hydropower constitutes a significant share of the electricity generation mix \citep{IEA2021Hydropower}.

Each hydropower plant in the cascade consists of hydraulic turbines for electricity generation, a diversion barrage for water routing, and a small-scale reservoir for short-term water storage management.
The stations are hydraulically coupled through the river network, such that upstream water releases directly affect downstream inflows and generation capabilities.

Wind and solar power outputs are modeled as exogenous stochastic processes, reflecting the inherent uncertainty and variability of meteorological conditions.
Consequently, vRES generation is non-dispatchable and may only be curtailed.
In contrast, the hydropower plants are fully dispatchable and represent the sole source of operational flexibility within the hybrid system.
We remark that, despite their limited storage capacity, run-of-the-river hydropower plants are widely recognized as particularly effective in mitigating vRES fluctuations while respecting hydraulic and operational constraints, as discussed in \cite{piron2016operating}.

In the following, the dispatch problem is formulated from the perspective of an energy producer operating the CH-vRES system, who is typically not directly responsible for grid management. 
Accordingly, network constraints are neglected in the present formulation and will be considered in future work.

The dispatch problem is formulated as a stochastic MPC problem, accounting for uncertainties in both water inflows within the cascaded hydropower system and vRES generation.
Forecasts of the uncertain parameters are denoted by the hat symbol, $\hat{\cdot}$,
and are represented using a finite scenario set $\boldsymbol{\Omega}$, indexed by $\omega \in \boldsymbol{\Omega}$,
with associated probabilities $\pi_{\omega}$, where $\sum_{\omega \in \boldsymbol{\Omega}} \pi_{\omega} = 1$.

Following the standard receding-horizon principle of MPC, updated forecasts are incorporated into the decision-making process at each time period $t \in \boldsymbol{T}$.
At time $t$, the controller solves an optimization problem over a prediction horizon $\boldsymbol{K}$, indexed by $k$, with cardinality $K \coloneqq |\boldsymbol{K}|$ and sampling time $\Delta$.
At the subsequent time period $t+1$, the prediction horizon $\boldsymbol{K}$ is shifted forward by one time period, and the optimization problem is resolved using updated forecasts.
A control decision computed at time $t$ for a future time period $t+k$ is denoted by $t+k \mid t$, representing the action planned for time $t+k$ based on the information available at time $t$.

\subsection{Full-Scale Centralized Stochastic Model Predictive Control Scheme}
\label{subsec:SMPC}

In this subsection, we formulate the energy dispatch problem as a conventional full-scale centralized stochastic MPC problem.
This formulation provides the foundation for the proposed temporal aggregation and scenario-based decomposition presented in the following subsections.

The CH-vRES hybrid system is modeled under the following assumptions:
(i) the tailrace water levels are assumed constant, owing to the negligible short-term influence of run-of-the-river plant discharges on downstream water elevations; and
(ii) all turbines within a given hydropower plant are assumed to be of the same type.
Under assumption (ii), the turbines are assumed to be homogeneous within each hydropower plant, while heterogeneous turbine characteristics are allowed across different plants in the cascade. 
This assumption is introduced primarily to simplify the notation and streamline the presentation of the proposed formulation.
Nevertheless, the adopted problem formulation can be naturally extended to account for heterogeneous turbines within the same plant, for instance by adopting the modeling approach proposed in \cite{anagnostopoulos2007optimal}.
We remark that the impact of these assumptions on the considered case study is evaluated in \cite{11270217}.

We consider a cascade of run-of-the-river hydropower plants indexed by $n \in \boldsymbol{N}$.
For plant $n$, the turbine and barrage discharges (m$^3$/s) at time $t+k$ in scenario $\omega$ are denoted by $q^{\mathrm{tr}}_{n,\omega,t+k|t}$ and $q^{\mathrm{br}}_{n,\omega,t+k|t}$, respectively,
while the reservoir forebay water level (m) is denoted by $l_{n,\omega,t+k|t}$.

The forebay water level dynamics of reservoir $n$ at time $t+k$ in scenario $\omega$ are governed by the reservoir surface area (m$^2$), denoted by $S_n$,
the water inflow (m$^3$/s) and outflow (m$^3$/s), denoted by $q^{\mathrm{in}}_{n,\omega,t+k|t}$ and $q^{\mathrm{out}}_{n,\omega,t+k|t}$, respectively,
and the current measured reservoir forebay water level (m) at time $t$, denoted by $L^0_{n}$, as follows:
\begin{align}
l_{n, \omega, t+k | t} & = l_{n, \omega, t+k-1 | t} + \frac{\left(q^{\mathrm{in}}_{n, \omega, t+k | t} - q^{\mathrm{out}}_{n, \omega, t+k | t}\right) \Delta}{S_n}, \; \forall n, \forall \omega, \forall k \in \boldsymbol{K} \setminus \{0\}, \label{eq:cen_sto_water_level}\\
l_{n, \omega, t | t} & = L^0_{n}, \; \forall n, \forall \omega, \label{eq:cen_sto_water_level_init}\\
l_{n, \omega, t + K-1 | t} & = L^0_{n}, \; \forall n, \forall \omega. \label{eq:cen_sto_water_level_final}
\end{align}

The inflow is given by the sum of upstream plant discharges and uncertain external inflows (m$^3$/s), denoted by $\hat{Q}^{\mathrm{ext}}_{n,\omega,t+k}$, from the river and its tributaries:
\begin{align}
    q^{\mathrm{in}}_{n, \omega, t+k | t} = & \; q^{\mathrm{br}}_{n-1, \omega, t+k | t} + q^{\mathrm{tr}}_{n-1, \omega, t+k | t} + \hat{Q}^\mathrm{ext}_{n, \omega, t+k}, \; \forall n \in \boldsymbol{N} \setminus \{0\}, \forall \omega, \forall k, \label{eq:cen_sto_inflow}\\
    q^{\mathrm{in}}_{0, \omega, t+k | t} = & \; \hat{Q}^\mathrm{ext}_{0, \omega, t+k}, \; \forall \omega, \forall k. \label{eq:cen_sto_inflow_init}
\end{align}
Similarly, the outflow is defined as:
\begin{equation}
\label{eq:cen_sto_outflow}
    q^{\mathrm{out}}_{n, \omega, t+k | t} = q^{\mathrm{br}}_{n, \omega, t+k | t} + q^{\mathrm{tr}}_{n, \omega, t+k | t}, \; \forall n, \forall \omega, \forall k.
\end{equation}

Water propagation dynamics can be incorporated into \eqref{eq:cen_sto_water_level}--\eqref{eq:cen_sto_outflow}.
Nevertheless, water travel times between consecutive hydropower plants are frequently neglected \cite{kumar2011genetic},
for instance because they are negligible relative to the problem sampling time \cite{xu2018analysis},
sufficiently detailed hydraulic data are unavailable \cite{li2026scheduling},
or the explicit modeling of such dynamics falls outside the scope of the study \cite{khodadadi2022stochastic}.
Other formulations represent water travel times as discrete parameters assumed to be known prior to solving the dispatch problem \cite{zheng2026leveraging,xie2021hydropower}.
Such an approximation may, however, be restrictive, as water travel times can vary continuously with hydraulic and operating conditions \cite{ge2018successive}.
More generally, the mismatch between the inflows computed while neglecting water travel times and the actual inflows observed at the hydropower plants can be represented as a stochastic disturbance within the controller \cite{zanoli2024model}.

Investigating alternative formulations for representing water propagation dynamics or assessing their impact on the hydropower dispatch decisions falls outside the scope of this work and is instead addressed in dedicated studies \cite{ge2018successive,yuan2022optimal};
hence, water travel times are not explicitly modeled in the present formulation.
This modeling choice allows the analysis to focus on a key methodological objective of this study:
handling nonconvex hydropower dynamics and intertemporal storage and ramping constraints within the proposed temporally aggregated distributed stochastic MPC scheme, while retaining a formal performance guarantee.
Assessing the impact of water travel times on the performance of the proposed controller is left for future work.

Importantly, this simplification does not preclude accounting for the residual modeling errors induced by unresolved water propagation dynamics.
Specifically, \eqref{eq:cen_sto_inflow} can be augmented with an additive disturbance capturing the mismatch between the inflows computed by the proposed controller and the actual inflows observed at the hydropower plants.
Such a disturbance can be modeled as a forecasted stochastic parameter, as in \cite{wu2025multi},
without affecting the methodological developments presented in this paper.
The proposed methodology therefore remains applicable when a disturbance model is employed to account for the residual effects induced by water travel times within the cascaded system.

To limit abrupt variations in the turbine operation and reduce mechanical stress and wear on both the turbines and the associated hydraulic infrastructure, a ramp limit (m$^3$/s), denoted by $\Delta^{\mathrm{tr}}_{n}$, is enforced on the turbine discharge:
\begin{align}
    q^{\mathrm{tr}}_{n, \omega, t+k|t} - q^{\mathrm{tr}}_{n, \omega, t+k-1|t} & \leq \Delta^{\mathrm{tr}}_{n}, \; \forall n, \forall \omega, \forall k \in \boldsymbol{K} \setminus \{0\}, \label{eq:cen_sto_ramping1}\\
    q^{\mathrm{tr}}_{n, \omega, t+k-1|t} - q^{\mathrm{tr}}_{n, \omega, t+k|t} & \leq \Delta^{\mathrm{tr}}_{n}, \; \forall n, \forall \omega, \forall k \in \boldsymbol{K} \setminus \{0\}. \label{eq:cen_sto_ramping2}
\end{align}

In addition to the ramping constraints in \eqref{eq:cen_sto_ramping1} and \eqref{eq:cen_sto_ramping2}, when activated, the turbines in the hydropower plant $n$ are subject to minimum and maximum discharge limits, denoted by $\underline{Q}^{\mathrm{tr}}_{n}$ and $\overline{Q}^{\mathrm{tr}}_{n}$, respectively.
Likewise, the hydropower generation of plant $n$ in scenario $\omega$ at time $t+k$, denoted by $p^{\mathrm{h}}_{n,\omega,t+k|t}$,
is bounded by minimum and maximum generation limits (MW), denoted by $\underline{P}^{\mathrm{h}}_{n}$ and $\overline{P}^{\mathrm{h}}_{n}$, respectively.
Thus, the operational limits on turbine discharges and hydropower generation are enforced as follows:
\begin{align}
    b_{n, \omega, t+k|t} \, \underline{Q}^{\mathrm{tr}}_{n} & \leq q^{\mathrm{tr}}_{n, \omega, t+k | t} \leq \overline{Q}^{\mathrm{tr}}_{n} \, b_{n, \omega, t+k|t}, \; \forall n, \forall \omega, \forall k, \label{eq:cen_sto_tr_lim}\\
    b_{n, \omega, t+k|t} \, \underline{P}^{\mathrm{h}}_{n} & \leq p^{\mathrm{h}}_{n, \omega, t+k | t} \leq \overline{P}^{\mathrm{h}}_{n} \, b_{n, \omega, t+k|t}, \; \forall n, \forall \omega, \forall k, \label{eq:cen_sto_ph_lim}\\
    b_{n, \omega, t+k|t} & \in \left\{0,1\right\}, \; \forall n, \forall \omega, \forall k. \label{eq:cen_sto_binary_lim}
\end{align}
Here, the binary variable $b_{n, \omega, t+k|t}$ indicates the on/off status of the turbines in hydropower plant $n$ in scenario $\omega$ at time $t+k$.

Moreover, let $\underline{L}_n$ and $\overline{L}_n$ denote the minimum and maximum forebay water levels (m) of reservoir $n$, respectively,
and let $\underline{Q}^{\mathrm{br}}_n$ denote the minimum barrage discharge (m$^3$/s).
The reservoir water levels and barrage discharges are subject to the following operational limits:
\begin{align}
    \underline{L}_{n} & \leq l_{n, \omega, t+k | t} \leq \overline{L}_{n}, \; \forall n, \forall \omega, \forall k, \label{eq:cen_sto_wl_lim}\\
    \underline{Q}^{\mathrm{br}}_{n} & \leq q^{\mathrm{br}}_{n, \omega, t+k | t}, \; \forall n, \forall \omega, \forall k. \label{eq:cen_sto_br_lim}
\end{align}

The power output of the $n$-th hydropower plant in scenario $\omega$ at time $t+k$ is a function of both the net hydraulic head (m), denoted by $h_{n,\omega,t+k|t}$, of the associated reservoir and the turbine discharge:
\begin{equation}
\label{eq:cen_sto_hydro_gen}
    p^{\mathrm{h}}_{n,\omega,t+k|t} = 10^{-6} \, w \, g \, \eta_n \, q^{\mathrm{tr}}_{n, \omega, t+k|t} \, h_{n, \omega, t+k|t}, \; \forall n, \forall \omega, \forall k,
\end{equation}
where $w$ denotes the water density (kg/m$^3$), $g$ denotes the gravitational acceleration (m/s$^2$), and $\eta_n$ denotes the turbine efficiency.
The factor $10^{-6}$ converts watts to megawatts.
The hydraulic head is defined as:
\begin{equation}
\label{eq:cen_sto_head}
    h_{n, \omega, t+k|t} = l_{n, \omega, t+k|t} - L^{\mathrm{tlr}}_{n}, \; \forall n, \forall \omega, \forall k,
\end{equation}
where $L^{\mathrm{tlr}}_n$ denotes the tailrace water level (m) of reservoir $n$.

The hydropower generation function \eqref{eq:cen_sto_hydro_gen} is inherently nonconvex due to its bilinear dependence on the turbine discharge and the hydraulic head.
Directly representing this bilinearity is often computationally intractable, particularly in the context of short-term, high-resolution dispatch strategies \citep{zhang2022chance}. 
As a result, the original formulation is commonly replaced by suitable approximations.
A variety of such approximation techniques exist in the literature,
the discussion of which lies beyond the scope of this paper and can be found in dedicated studies \citep{santos2022piecewise}.

In this paper, we adopt the widely used McCormick relaxation \citep{zhang2022long}, whereby the bilinear term in \eqref{eq:cen_sto_hydro_gen} is replaced by a convex envelope.
Given that the minimum turbine discharge permitted by the operational constraints \eqref{eq:cen_sto_tr_lim} is zero (corresponding to the turbine being shut down),
the McCormick relaxation of \eqref{eq:cen_sto_hydro_gen} simplifies to the following form:
\begin{align}
p^{\mathrm{h}}_{n, \omega, t+k|t} & \geq 10^{-6} \, w \, g \, \eta_n \, \underline{H}_{n} \, q^{\mathrm{tr}}_{n, \omega, t+k|t}, \; \forall n, \forall \omega, \forall k, \label{eq:cen_sto_McCormick1}\\
p^{\mathrm{h}}_{n, \omega, t+k|t} & \geq 10^{-6} \, w \, g \, \eta_n \left(\overline{Q}^{\mathrm{tr}}_{n} \, h_{n, \omega, t+k|t} + \overline{H}_{n} \, q^{\mathrm{tr}}_{n, \omega, t+k|t} - \overline{Q}^{\mathrm{tr}}_{n}  \, \overline{H}_{n}\right), \; \forall n, \forall \omega, \forall k, \label{eq:cen_sto_McCormick2}\\
p^{\mathrm{h}}_{n, \omega, t+k|t} & \leq 10^{-6} \, w \, g \, \eta_n \, \overline{H}_{n} \, q^{\mathrm{tr}}_{n, \omega, t+k|t}, \; \forall n, \forall \omega, \forall k, \label{eq:cen_sto_McCormick3}\\
p^{\mathrm{h}}_{n, \omega, t+k|t} & \leq 10^{-6} \, w \, g \, \eta_n \left(\overline{Q}^{\mathrm{tr}}_{n} \, h_{n, \omega, t+k|t} + \underline{H}_{n} \, q^{\mathrm{tr}}_{n, \omega, t+k|t} - \overline{Q}^{\mathrm{tr}}_n \, \underline{H}_{n}\right), \; \forall n, \forall \omega, \forall k. \label{eq:cen_sto_McCormick4}
\end{align}
Here,
$\underline{H}_n \coloneqq \underline{L}_n - L^{\mathrm{tlr}}_n$
and
$\overline{H}_n \coloneqq \overline{L}_n - L^{\mathrm{tlr}}_n$
denote the minimum and maximum head values (m), respectively, of hydropower plant $n$.

The McCormick relaxation is among the most widely adopted relaxation techniques for bilinear functions in general \citep{najman2019tightness}, and for hydropower modeling in particular \citep{shi2024stochastic}.
It is widely regarded as one of the most accurate convex relaxations for bilinear terms \citep{blom2024single}.
A theoretical analysis of this relaxation in the context of hydropower modeling is presented in \cite{flamm2020two}, while a numerical assessment for the case study considered in this work is provided in \cite{11270217}.

Wind and solar power generation (MW) in scenario $\omega$ at time $t+k$ are denoted by $p^{\mathrm{w}}_{\omega,t+k|t}$ and $p^{\mathrm{s}}_{\omega,t+k|t}$, respectively,
and are bounded by their uncertain capacity factors $\hat{F}^{\mathrm{W}}_{\omega,t+k}$ and $\hat{F}^{\mathrm{S}}_{\omega,t+k}$.
This yields the following constraints:
\begin{align}
    0 & \leq p^{\mathrm{w}}_{\omega, t+k|t} \leq \hat{F}^{\mathrm{W}}_{\omega, t+k} X^{\mathrm{W}}, \; \forall \omega, \forall k, \label{eq:cen_sto_wind}\\
    0 & \leq p^{\mathrm{s}}_{\omega, t+k|t} \leq \hat{F}^{\mathrm{S}}_{\omega, t+k} X^{\mathrm{S}}, \; \forall \omega, \forall k, \label{eq:cen_sto_solar}
\end{align}
where $X^{\mathrm{W}}$ and $X^{\mathrm{S}}$ denote the installed capacities (MW) of the wind and solar generation units, respectively.

Then, the total power output (MW) of the CH-vRES hybrid system in scenario $\omega$ at time $t+k$, denoted by $p_{\omega,t+k|t}$, is defined as the sum of hydropower, wind, and solar generation:
\begin{equation}
    p_{\omega, t+k|t} = \sum_{n \in \boldsymbol{N}} p^{\mathrm{h}}_{n, \omega, t+k|t} + p^{\mathrm{w}}_{\omega, t+k|t} + p^{\mathrm{s}}_{\omega, t+k|t}, \; \forall \omega, \forall k. \label{eq:cen_sto_power_balance}
\end{equation}

Let $\boldsymbol{u}^{\boldsymbol{\mathrm{h}}}_{n,t|t}$ and $\boldsymbol{u}^{\boldsymbol{\mathrm{vRES}}}_{t|t}$ denote the vectors of control actions for the hydropower plant $n$ and the vRES, respectively, computed by the MPC scheme at time $t$ for $k=0$.
The decision variables of the dispatch model are collected in the set $\boldsymbol{z}$, defined as:
\begin{equation}
\begin{aligned}
\boldsymbol{z} \coloneqq \Big\{
& \boldsymbol{u}^{\boldsymbol{\mathrm{h}}}_{n,t|t},
\boldsymbol{u}^{\boldsymbol{\mathrm{vRES}}}_{t|t},
p_{\omega,t+k|t},
p^{\mathrm{h}}_{n,\omega,t+k|t},
p^{\mathrm{w}}_{\omega,t+k|t},
p^{\mathrm{s}}_{\omega,t+k|t},
l_{n,\omega,t+k|t},
q^{\mathrm{in}}_{n,\omega,t+k|t}, \\
& q^{\mathrm{out}}_{n,\omega,t+k|t},
q^{\mathrm{tr}}_{n,\omega,t+k|t},
q^{\mathrm{br}}_{n,\omega,t+k|t},
h_{n,\omega,t+k|t},
b_{n,\omega,t+k|t}
\;\Big|\;
n \in \boldsymbol{N},
\omega \in \boldsymbol{\Omega},
k \in \boldsymbol{K}
\Big\}.
\end{aligned}
\end{equation}

The objective of the dispatch problem is to track a reference power signal (MW), denoted by $P^{\mathrm{ref}}_{t+k}$.
To this end, the \textbf{full-scale centralized stochastic MPC scheme} solves, at time $t$ over the prediction horizon $\boldsymbol{K}$, the following MIQP problem:
\begin{subequations}
\label{eq:cen_sto_MPC}
\begin{align}
\min_{\boldsymbol{z}} \quad & F(\boldsymbol{z}) \coloneqq \sum_{\omega \in \boldsymbol{\Omega}} \pi_{\omega} \sum_{k \in \boldsymbol{K}} \left( p_{\omega, t+k|t} - P^{\mathrm{ref}}_{t+k} \right)^2 \label{eq:cen_sto_obj}\\
\text{s.t.} \quad & \eqref{eq:cen_sto_water_level}-\eqref{eq:cen_sto_br_lim}, \eqref{eq:cen_sto_head}-\eqref{eq:cen_sto_power_balance}\\
& \left[
q^{\mathrm{br}}_{n,\omega,t|t},
q^{\mathrm{tr}}_{n,\omega,t|t}
\right]^\top = \boldsymbol{u}^{\boldsymbol{\mathrm{h}}}_{n, t|t}, \; \forall n, \forall \omega, \label{eq:cen_sto_MPC_control_hydro}\\
& \left[
p^{\mathrm{w}}_{\omega, t|t},
p^{\mathrm{s}}_{\omega, t|t}
\right]^\top = \boldsymbol{u}^{\boldsymbol{\mathrm{vRES}}}_{t|t}, \; \forall \omega. \label{eq:cen_sto_MPC_control_vres}
\end{align}
\end{subequations}
Following the rolling-horizon approach, only the first control actions,
$\boldsymbol{u}^{\boldsymbol{\mathrm{h}}}_{n,t|t}$ for all $n$ and $\boldsymbol{u}^{\boldsymbol{\mathrm{vRES}}}_{t|t}$,
are implemented at each time period, and the prediction horizon $\boldsymbol{K}$ is shifted forward by one time period at time $t+1$.

We remark that reference-tracking terms of the form \eqref{eq:cen_sto_obj} are standard in MPC formulations \citep{bordons2019model} and provide a flexible modeling framework for representing a wide range of operational objectives relevant to energy producers.
In particular, such formulations can capture the real-time tracking of bidding curves previously submitted to electricity markets, with the aim of reducing imbalance penalties or providing ancillary services \citep{santosuosso2024stochastic},
as well as the tracking of reference signals issued by system operators for grid management purposes, such as peak shaving \citep{dongol2018model}.

\subsection{Temporally Aggregated Centralized Stochastic Model Predictive Control Scheme}
\label{subsec:SMPC_TSA}

Solving the nonconvex stochastic energy dispatch model \eqref{eq:cen_sto_MPC} at high temporal resolution within real-time dispatch strategies is often computationally intractable.
To alleviate this complexity, TSA is employed to construct a temporally aggregated counterpart of the full-scale model \eqref{eq:cen_sto_MPC},
defined over a reduced set of representative time periods (or clusters) $\boldsymbol{R}$, indexed by $r$, with cardinality $R \coloneqq |\boldsymbol{R}|$.
When the number of clusters satisfies $R \ll K$, the temporally aggregated model yields significant computational savings.

We group the decision variables of the temporally aggregated model in $\boldsymbol{\bar{z}}$, defined as
\begin{equation}
\begin{aligned}
\boldsymbol{\bar{z}} \coloneqq \Big\{
& \boldsymbol{\bar{u}}^{\boldsymbol{\mathrm{h}}}_{n, t|t},
\boldsymbol{\bar{u}}^{\boldsymbol{\mathrm{vRES}}}_{t|t},
\bar{p}_{\omega,t+r|t},
\bar{p}^{\mathrm{h}}_{n,\omega,t+r|t},
\bar{p}^{\mathrm{w}}_{\omega,t+r|t},
\bar{p}^{\mathrm{s}}_{\omega,t+r|t},
\bar{l}_{n,\omega,t+r|t},
\bar{q}^{\mathrm{in}}_{n,\omega,t+r|t}, \\
& \bar{q}^{\mathrm{out}}_{n,\omega,t+r|t},
\bar{q}^{\mathrm{tr}}_{n,\omega,t+r|t},
\bar{q}^{\mathrm{br}}_{n,\omega,t+r|t},
\bar{h}_{n,\omega,t+r|t},
\bar{b}_{n,\omega,t+r|t}
\;\Big|\;
n \in \boldsymbol{N},
\omega \in \boldsymbol{\Omega},
r \in \boldsymbol{R}
\Big\}.
\end{aligned}
\end{equation}

Let $\boldsymbol{K}_r \subseteq \boldsymbol{K}$ denote the set of time periods $k \in \boldsymbol{K}$ assigned to cluster $r \in \boldsymbol{R}$,
with cardinality $K_r \coloneqq |\boldsymbol{K}_r|$.
Then, the \textbf{temporally aggregated centralized stochastic MPC scheme} solves the following MIQP problem at time $t$, defined over the set of representative time periods $\boldsymbol{R}$:
\begin{subequations}
\label{eq:TSA_cen_sto_MPC}
\begin{align}
\min_{\boldsymbol{\bar{z}}} \quad &
\bar{F}(\boldsymbol{\bar{z}}) \coloneqq
\sum_{\omega \in \boldsymbol{\Omega}} \pi_{\omega} \sum_{r \in \boldsymbol{R}} K_r \, \left( \bar{p}_{\omega, t+r|t} - \frac{1}{K_r} \sum_{k \in \boldsymbol{K}_r} P^{\mathrm{ref}}_{t+k} \right)^2 \label{eq:TSA_cen_sto_obj}\\
\text{s.t.} \quad &
\bar{l}_{n, \omega, t+r | t} = \bar{l}_{n, \omega, t+r-1 | t} + K_r \; \Delta \; \frac{ \bar{q}^{\mathrm{in}}_{n, \omega, t+r | t} - \bar{q}^{\mathrm{out}}_{n, \omega, t+r | t}}{S_n}, \; \forall n, \forall \omega, \forall r \in \boldsymbol{R} \setminus \{0\},
\label{eq:TSA_cen_sto_water_level}\\
& \bar{l}_{n, \omega, t | t} = L^0_{n}, \; \forall n, \forall \omega,
\label{eq:TSA_cen_sto_water_level_init}\\
& \bar{l}_{n, \omega, t + R-1 | t} = L^0_{n}, \; \forall n, \forall \omega, \label{eq:TSA_cen_sto_water_level_final}\\
& \bar{q}^{\mathrm{in}}_{n, \omega, t+r | t} = \bar{q}^{\mathrm{br}}_{n-1, \omega, t+r | t} + \bar{q}^{\mathrm{tr}}_{n-1, \omega, t+r | t} + \frac{1}{K_r} \sum_{k \in \boldsymbol{K}_r} \hat{Q}^\mathrm{ext}_{n, \omega, t+k}, \; \forall n \in \boldsymbol{N} \setminus \{0\}, \forall \omega, \forall r,
\label{eq:TSA_cen_sto_inflow}\\
& \bar{q}^{\mathrm{in}}_{0, \omega, t+r | t} = \frac{1}{K_r} \sum_{k \in \boldsymbol{K}_r} \hat{Q}^\mathrm{ext}_{0, \omega, t+k}, \; \forall \omega, \forall r,
\label{eq:TSA_cen_sto_inflow_init}\\
& \bar{q}^{\mathrm{out}}_{n, \omega, t+r | t} = \bar{q}^{\mathrm{br}}_{n, \omega, t+r | t} + \bar{q}^{\mathrm{tr}}_{n, \omega, t+r | t}, \; \forall n, \forall \omega, \forall r,
\label{eq:TSA_cen_sto_outflow}\\
& \bar{q}^{\mathrm{tr}}_{n, \omega, t+r|t} - \bar{q}^{\mathrm{tr}}_{n, \omega, t+r-1|t} \, K_{r-1} \leq
\Delta^{\mathrm{tr}}_{n} + \Delta^{\mathrm{tr}}_{n} \; \frac{K_r - 1}{2}, \; \forall n, \forall \omega, \forall r \in \boldsymbol{R} \setminus \{0\}, 
\label{eq:TSA_cen_sto_ramping1}\\
& \bar{q}^{\mathrm{tr}}_{n, \omega, t+r-1|t} - \bar{q}^{\mathrm{tr}}_{n, \omega, t+r|t} \, K_{r} \leq
\Delta^{\mathrm{tr}}_{n} + \Delta^{\mathrm{tr}}_{n} \; \frac{K_{r - 1} - 1}{2}, \; \forall n, \forall \omega, \forall r \in \boldsymbol{R} \setminus \{0\}, \label{eq:TSA_cen_sto_ramping2}\\
& \bar{b}_{n, \omega, t+r|t} \, \underline{Q}^{\mathrm{tr}}_{n} \leq \bar{q}^{\mathrm{tr}}_{n, \omega, t+r | t} \leq \overline{Q}^{\mathrm{tr}}_{n} \, \bar{b}_{n, \omega, t+r|t}, \; \forall n, \forall \omega, \forall r, 
\label{eq:TSA_cen_sto_tr_lim}\\
& \bar{b}_{n, \omega, t+r|t} \, \underline{P}^{\mathrm{h}}_{n} \leq \bar{p}^{\mathrm{h}}_{n, \omega, t+r | t} \leq \overline{P}^{\mathrm{h}}_{n} \, \bar{b}_{n, \omega, t+r|t}, \; \forall n, \forall \omega, \forall r, 
\label{eq:TSA_cen_sto_ph_lim}\\
& \underline{L}_{n} \leq \bar{l}_{n, \omega, t+r | t} \leq \overline{L}_{n}, \; \forall n, \forall \omega, \forall r,
\label{eq:TSA_cen_sto_wl_lim}\\
& \underline{Q}^{\mathrm{br}}_{n} \leq \bar{q}^{\mathrm{br}}_{n, \omega, t+r | t}, \; \forall n, \forall \omega, \forall r,
\label{eq:TSA_cen_sto_br_lim}\\
& \bar{h}_{n, \omega, t+r|t} = \bar{l}_{n, \omega, t+r|t} - L^{\mathrm{tlr}}_{n},
\; \forall n, \forall \omega, \forall r,
\label{eq:TSA_cen_sto_head}\\
& \bar{p}^{\mathrm{h}}_{n, \omega, t+r|t} \geq 10^{-6} \, w \, g \, \eta_n \, \underline{H}_{n} \, \bar{q}^{\mathrm{tr}}_{n, \omega, t+r|t}, \; \forall n, \forall \omega, \forall r, 
\label{eq:TSA_cen_sto_McCormick1}\\
& \bar{p}^{\mathrm{h}}_{n, \omega, t+r|t} \geq 10^{-6} \, w \, g \, \eta_n \left(
\overline{Q}^{\mathrm{tr}}_{n}
\left(\frac{\bar{h}_{n, \omega, t+r|t}}{K_r} + \frac{K_r - 1}{K_r} \, \underline{H}_n \right)
+ \overline{H}_{n} \, \bar{q}^{\mathrm{tr}}_{n, \omega, t+r|t} - \overline{Q}^{\mathrm{tr}}_{n}  \, \overline{H}_{n}\right), \; \forall n, \forall \omega, \forall r, 
\label{eq:TSA_cen_sto_McCormick2}\\
& \bar{p}^{\mathrm{h}}_{n, \omega, t+r|t} \leq 10^{-6} \, w \, g \, \eta_n \, \overline{H}_{n} \, \bar{q}^{\mathrm{tr}}_{n, \omega, t+r|t}, \; \forall n, \forall \omega, \forall r,
\label{eq:TSA_cen_sto_McCormick3}\\
& \bar{p}^{\mathrm{h}}_{n, \omega, t+r|t} \leq 10^{-6} \, w \, g \, \eta_n \left(\overline{Q}^{\mathrm{tr}}_{n} 
\left( \frac{\bar{h}_{n, \omega, t+r|t}}{K_r} + \frac{K_r - 1}{K_r} \, \overline{H}_n \right)
+ \underline{H}_{n} \, \bar{q}^{\mathrm{tr}}_{n, \omega, t+r|t} - \overline{Q}^{\mathrm{tr}}_n \, \underline{H}_{n}\right), \; \forall n, \forall \omega, \forall r, 
\label{eq:TSA_cen_sto_McCormick4}\\
& 0 \leq \bar{p}^{\mathrm{w}}_{\omega, t+r|t} \leq \frac{X^{\mathrm{W}}}{K_r} \, \sum_{k \in \boldsymbol{K}_r} \hat{F}^{\mathrm{W}}_{\omega, t+k}, \; \forall \omega, \forall r,
\label{eq:TSA_cen_sto_wind}\\
& 0 \leq \bar{p}^{\mathrm{s}}_{\omega, t+r|t} \leq \frac{X^{\mathrm{S}}}{K_r} \sum_{k \in \boldsymbol{K}_r} \hat{F}^{\mathrm{S}}_{\omega, t+k}, \; \forall \omega, \forall r,
\label{eq:TSA_cen_sto_solar}\\
& \bar{p}_{\omega, t+r|t} = \sum_{n \in \boldsymbol{N}} \bar{p}^{\mathrm{h}}_{n, \omega, t+r|t} + \bar{p}^{\mathrm{w}}_{\omega, t+r|t} + \bar{p}^{\mathrm{s}}_{\omega, t+r|t}, \; \forall \omega, \forall r, 
\label{eq:TSA_cen_sto_power_balance}\\
& \left[
\bar{q}^{\mathrm{br}}_{n,\omega,t|t},
\bar{q}^{\mathrm{tr}}_{n,\omega,t|t}
\right]^\top = \boldsymbol{\bar{u}}^{\boldsymbol{\mathrm{h}}}_{n, t|t}, \; \forall n, \forall \omega,
\label{eq:TSA_cen_sto_MPC_control_hydro}\\
& \left[
\bar{p}^{\mathrm{w}}_{\omega, t|t},
\bar{p}^{\mathrm{s}}_{\omega, t|t}
\right]^\top = \boldsymbol{\bar{u}}^{\boldsymbol{\mathrm{vRES}}}_{t|t}, \; \forall \omega,
\label{eq:TSA_cen_sto_MPC_control_vres}\\
& \bar{b}_{n,\omega,t+r\mid t} \in \left\{ \frac{m}{K_r} \,\middle|\, m = 0,1,\dots,K_r \right\}, \; \forall n, \forall \omega, \forall r.
\label{eq:TSA_cen_sto_binary_lim}
\end{align}
\end{subequations}

In \eqref{eq:TSA_cen_sto_MPC},
\eqref{eq:TSA_cen_sto_obj} denotes the aggregated counterpart of the original full-scale objective function \eqref{eq:cen_sto_obj},
while the constraints \eqref{eq:TSA_cen_sto_water_level}--\eqref{eq:TSA_cen_sto_binary_lim} represents the aggregated counterparts of the full-scale constraints
\eqref{eq:cen_sto_water_level}--\eqref{eq:cen_sto_ph_lim}, \eqref{eq:cen_sto_wl_lim}, \eqref{eq:cen_sto_br_lim}, \eqref{eq:cen_sto_head}--\eqref{eq:cen_sto_power_balance}, \eqref{eq:cen_sto_MPC_control_hydro}, \eqref{eq:cen_sto_MPC_control_vres}, and \eqref{eq:cen_sto_binary_lim}, respectively.

Notably, most aggregated constraints preserve the structural form of their full-scale counterparts,
differing only in that they are defined over the
set of representative time periods $\boldsymbol{R}$ rather than $\boldsymbol{K}$.
In contrast, the intertemporal and nonconvex constraints require a dedicated treatment.
Specifically, the aggregated storage constraints \eqref{eq:TSA_cen_sto_water_level} and the ramping constraints \eqref{eq:TSA_cen_sto_ramping1}--\eqref{eq:TSA_cen_sto_ramping2} are augmented,
relative to their full-scale counterparts in \eqref{eq:cen_sto_water_level}, and \eqref{eq:cen_sto_ramping1}--\eqref{eq:cen_sto_ramping2}, respectively,
with cluster-specific terms.
Similarly, the aggregated constraints \eqref{eq:TSA_cen_sto_McCormick2} and \eqref{eq:TSA_cen_sto_McCormick4},
which are directly affected by the intertemporal constraints \eqref{eq:TSA_cen_sto_water_level} through \eqref{eq:TSA_cen_sto_head},
are reformulated relative to their full-scale counterparts \eqref{eq:cen_sto_McCormick2} and \eqref{eq:cen_sto_McCormick4} by incorporating additional cluster-specific terms.
Moreover, the binary decision variables $b_{n,\omega,t+k\mid t}$ of the full-scale model are reformulated in the aggregated model as $\bar{b}_{n,\omega,t+r\mid t}$,
defined as variables taking values in a finite, cluster-specific discrete set, as enforced by \eqref{eq:TSA_cen_sto_binary_lim}.
These reformulations are introduced to establish the theoretical result presented in the following proposition.

\begin{proposition}\label{prop:main_result}
Let $\boldsymbol{z}$ be a feasible solution of the full-scale centralized stochastic model \eqref{eq:cen_sto_MPC}.
Let $\boldsymbol{\bar{z}}$ be derived from $\boldsymbol{z}$ as follows:

\noindent
\begin{minipage}{0.5\textwidth}
\begin{equation}\label{prop_agg_vars1}
\bar{p}_{\omega,t+r|t} \coloneqq
\frac{1}{K_r} \sum_{k \in \boldsymbol{K}_r} p_{\omega,t+k|t},
\; \forall \omega, \forall r,
\end{equation}
\end{minipage}%
\hfill
\begin{minipage}{0.5\textwidth}
\begin{equation}\label{prop_agg_vars2}
\bar{p}^{\mathrm{h}}_{n,\omega,t+r|t} \coloneqq
\frac{1}{K_r} \sum_{k \in \boldsymbol{K}_r}
p^{\mathrm{h}}_{n,\omega,t+k|t},
\; \forall n, \forall \omega, \forall r,
\end{equation}
\end{minipage}

\noindent
\begin{minipage}{0.5\textwidth}
\begin{equation}\label{prop_agg_vars3}
\bar{p}^{\mathrm{w}}_{\omega,t+r|t} \coloneqq
\frac{1}{K_r} \sum_{k \in \boldsymbol{K}_r}
p^{\mathrm{w}}_{\omega,t+k|t},
\; \forall \omega, \forall r,
\end{equation}
\end{minipage}%
\hfill
\begin{minipage}{0.5\textwidth}
\begin{equation}\label{prop_agg_vars4}
\bar{p}^{\mathrm{s}}_{\omega,t+r|t} \coloneqq
\frac{1}{K_r} \sum_{k \in \boldsymbol{K}_r}
p^{\mathrm{s}}_{\omega,t+k|t},
\; \forall \omega, \forall r,
\end{equation}
\end{minipage}

\noindent
\begin{minipage}{0.5\textwidth}
\begin{equation}\label{prop_agg_vars5}
\bar{q}^{\mathrm{in}}_{n,\omega,t+r|t} \coloneqq
\frac{1}{K_r} \sum_{k \in \boldsymbol{K}_r}
q^{\mathrm{in}}_{n,\omega,t+k|t},
\; \forall n, \forall \omega, \forall r,
\end{equation}
\end{minipage}%
\hfill
\begin{minipage}{0.5\textwidth}
\begin{equation}\label{prop_agg_vars6}
\bar{q}^{\mathrm{out}}_{n,\omega,t+r|t} \coloneqq
\frac{1}{K_r} \sum_{k \in \boldsymbol{K}_r}
q^{\mathrm{out}}_{n,\omega,t+k|t},
\; \forall n, \forall \omega, \forall r,
\end{equation}
\end{minipage}

\noindent
\begin{minipage}{0.5\textwidth}
\begin{equation}\label{prop_agg_vars7}
\bar{q}^{\mathrm{tr}}_{n,\omega,t+r|t} \coloneqq
\frac{1}{K_r} \sum_{k \in \boldsymbol{K}_r}
q^{\mathrm{tr}}_{n,\omega,t+k|t},
\; \forall n, \forall \omega, \forall r,
\end{equation}
\end{minipage}%
\hfill
\begin{minipage}{0.5\textwidth}
\begin{equation}\label{prop_agg_vars8}
\bar{q}^{\mathrm{br}}_{n,\omega,t+r|t} \coloneqq
\frac{1}{K_r} \sum_{k \in \boldsymbol{K}_r}
q^{\mathrm{br}}_{n,\omega,t+k|t},
\; \forall n, \forall \omega, \forall r,
\end{equation}
\end{minipage}

\noindent
\begin{minipage}{0.5\textwidth}
\begin{equation}\label{prop_agg_vars9}
\bar{h}_{n,\omega,t+r|t} \coloneqq
h_{n,\omega,t+\min_{k \in \boldsymbol{K}_r}(t+k)|t},
\; \forall n, \forall \omega, \forall r,
\end{equation}
\end{minipage}%
\hfill
\begin{minipage}{0.5\textwidth}
\begin{equation}\label{prop_agg_vars10}
\bar{b}_{n,\omega,t+r|t} \coloneqq
\frac{1}{K_r} \sum_{k \in \boldsymbol{K}_r}
b_{n,\omega,t+k|t},
\; \forall n, \forall \omega, \forall r,
\end{equation}
\end{minipage}

\noindent
\begin{minipage}{0.5\textwidth}
\begin{equation}\label{prop_agg_vars11}
\bar{l}_{n,\omega,t+r|t} \coloneqq
l_{n,\omega, \min_{k \in \boldsymbol{K}_r}(t+k)|t},
\; \forall n, \forall \omega, \forall r.
\end{equation}
\end{minipage}%
\hfill
\begin{minipage}{0.5\textwidth}
\end{minipage}
\vspace{0.5cm}

Moreover, let the temporally aggregated centralized stochastic dispatch model \eqref{eq:TSA_cen_sto_MPC} be constructed by clustering the original set of time periods $\boldsymbol{K}$ under the following assumptions:
\begin{itemize}
\item[A1:] The clustering preserves temporal chronology,
i.e., for any $k' \in \boldsymbol{K}_{r'}$ and $k'' \in \boldsymbol{K}_{r''}$ with $r' < r''$, it holds that $k' < k''$.
\item[A2:] The first and last time periods of the prediction horizon $\boldsymbol{K}$ are preserved as singleton clusters, namely,
$|\boldsymbol{K}_0| = 1$ and
$|\boldsymbol{K}_{R-1}| = 1$.
\end{itemize}

Then, $\boldsymbol{\bar{z}}$ is a feasible solution for the temporally aggregated centralized stochastic dispatch model
\eqref{eq:TSA_cen_sto_MPC}, and it holds that $\bar{F}\left(\boldsymbol{\bar{z}}\right) \leq F\left(\boldsymbol{z}\right)$.
\end{proposition}

\begin{proof}
We first demonstrate that any $\boldsymbol{\bar{z}}$ obtained through \eqref{prop_agg_vars1}--\eqref{prop_agg_vars11} is a feasible solution for the temporally aggregated model \eqref{eq:TSA_cen_sto_MPC}.

From \eqref{prop_agg_vars5}, \eqref{prop_agg_vars7} and \eqref{prop_agg_vars8},
substituting $\boldsymbol{\bar{z}}$ into
\eqref{eq:TSA_cen_sto_inflow}, \eqref{eq:TSA_cen_sto_inflow_init} and \eqref{eq:TSA_cen_sto_outflow}
of the temporally aggregated model \eqref{eq:TSA_cen_sto_MPC} yields
\begin{align}
\sum_{k \in \boldsymbol{K}_r} q^{\mathrm{in}}_{n,\omega,t+k|t} & =
\sum_{k \in \boldsymbol{K}_r} q^{\mathrm{br}}_{n-1,\omega,t+k|t}
+
\sum_{k \in \boldsymbol{K}_r} q^{\mathrm{tr}}_{n-1,\omega,t+k|t}
+ \sum_{k \in \boldsymbol{K}_r} \hat{Q}^\mathrm{ext}_{n, \omega, t+k},
\; \forall n \in \boldsymbol{N} \setminus \{0\}, \forall \omega, \forall r,
\label{prop:TSA_cen_sto_inflow}\\
\sum_{k \in \boldsymbol{K}_r} q^{\mathrm{in}}_{0,\omega,t+k|t} & =
\sum_{k \in \boldsymbol{K}_r} \hat{Q}^\mathrm{ext}_{0, \omega, t+k},
\; \forall \omega, \forall r,
\label{prop:TSA_cen_sto_inflow_init}\\
\sum_{k \in \boldsymbol{K}_r} q^{\mathrm{out}}_{n,\omega,t+k|t} & =
\sum_{k \in \boldsymbol{K}_r} q^{\mathrm{br}}_{n,\omega,t+k|t}
+
\sum_{k \in \boldsymbol{K}_r} q^{\mathrm{tr}}_{n,\omega,t+k|t},
\; \forall n, \forall \omega, \forall r.
\label{prop:TSA_cen_sto_outflow}
\end{align}

From \eqref{prop_agg_vars2}, \eqref{prop_agg_vars7}, \eqref{prop_agg_vars8} and \eqref{prop_agg_vars10},
substituting $\boldsymbol{\bar{z}}$ into
\eqref{eq:TSA_cen_sto_tr_lim}, \eqref{eq:TSA_cen_sto_ph_lim} and \eqref{eq:TSA_cen_sto_br_lim}
of the temporally aggregated model \eqref{eq:TSA_cen_sto_MPC} yields
\begin{align}
\underline{Q}^{\mathrm{tr}}_{n}
\sum_{k \in \boldsymbol{K}_r} b_{n,\omega,t+k|t}
& \leq
\sum_{k \in \boldsymbol{K}_r} q^{\mathrm{tr}}_{n,\omega,t+k|t}
\leq
\overline{Q}^{\mathrm{tr}}_{n}
\sum_{k \in \boldsymbol{K}_r} b_{n,\omega,t+k|t},
\; \forall n, \forall \omega, \forall r,
\label{prop:TSA_cen_sto_tr_lim}\\
\underline{P}^{\mathrm{h}}_{n}
\sum_{k \in \boldsymbol{K}_r} b_{n,\omega,t+k|t}
& \leq
\sum_{k \in \boldsymbol{K}_r} p^{\mathrm{h}}_{n,\omega,t+k|t}
\leq
\overline{P}^{\mathrm{h}}_{n}
\sum_{k \in \boldsymbol{K}_r} b_{n,\omega,t+k|t},
\; \forall n, \forall \omega, \forall r,
\label{prop:TSA_cen_sto_ph_lim}\\
K_r \, \underline{Q}^{\mathrm{br}}_{n} & \leq
\sum_{k \in \boldsymbol{K}_r}
q^{\mathrm{br}}_{n, \omega, t+k | t},
\; \forall n, \forall \omega, \forall r.
\label{prop:TSA_cen_sto_br_lim}
\end{align}

From \eqref{prop_agg_vars3} and \eqref{prop_agg_vars4},
substituting $\boldsymbol{\bar{z}}$ into
\eqref{eq:TSA_cen_sto_wind} and \eqref{eq:TSA_cen_sto_solar}
of the temporally aggregated model \eqref{eq:TSA_cen_sto_MPC} yields
\begin{align}
0 & \leq
\sum_{k \in \boldsymbol{K}_r} p^{\mathrm{w}}_{\omega, t+k|t}
\leq X^{\mathrm{W}} \sum_{k \in \boldsymbol{K}_r} \hat{F}^{\mathrm{W}}_{\omega, t+k},
\; \forall \omega, \forall r,
\label{prop:TSA_cen_sto_wind}\\
0 & \leq
\sum_{k \in \boldsymbol{K}_r} p^{\mathrm{s}}_{\omega, t+k|t}
\leq X^{\mathrm{S}} \sum_{k \in \boldsymbol{K}_r} \hat{F}^{\mathrm{S}}_{\omega, t+k},
\; \forall \omega, \forall r.
\label{prop:TSA_cen_sto_solar}
\end{align}

From \eqref{prop_agg_vars1}--\eqref{prop_agg_vars4},
substituting $\boldsymbol{\bar{z}}$ into
\eqref{eq:TSA_cen_sto_power_balance}
of the temporally aggregated model \eqref{eq:TSA_cen_sto_MPC} yields
\begin{equation}
\sum_{k \in \boldsymbol{K}_r} p_{\omega, t+k|t}
=
\sum_{k \in \boldsymbol{K}_r} \left(
\sum_{n \in \boldsymbol{N}} p^{\mathrm{h}}_{n, \omega, t+k|t} + p^{\mathrm{w}}_{\omega, t+k|t} + p^{\mathrm{s}}_{\omega, t+k|t}
\right), \; \forall \omega, \forall r.
\label{prop:TSA_cen_sto_power_balance}
\end{equation}

The full-scale constraints
\eqref{eq:cen_sto_inflow}--\eqref{eq:cen_sto_outflow},
\eqref{eq:cen_sto_tr_lim}, \eqref{eq:cen_sto_ph_lim}, \eqref{eq:cen_sto_br_lim},
\eqref{eq:cen_sto_wind}, \eqref{eq:cen_sto_solar}, and \eqref{eq:cen_sto_power_balance},
which are imposed individually for each $k \in \boldsymbol{K}$,
directly imply the corresponding aggregated constraints
\eqref{prop:TSA_cen_sto_inflow}--\eqref{prop:TSA_cen_sto_outflow},
\eqref{prop:TSA_cen_sto_tr_lim}, \eqref{prop:TSA_cen_sto_ph_lim}, \eqref{prop:TSA_cen_sto_br_lim},
\eqref{prop:TSA_cen_sto_wind}, \eqref{prop:TSA_cen_sto_solar}, and \eqref{prop:TSA_cen_sto_power_balance},
which are obtained by summing the corresponding full-scale constraints over all time periods $k$ assigned to each cluster $\boldsymbol{K}_r$.

From \eqref{prop_agg_vars2}, \eqref{prop_agg_vars7} and \eqref{prop_agg_vars9},
substituting $\boldsymbol{\bar{z}}$ into
\eqref{eq:TSA_cen_sto_McCormick1}--\eqref{eq:TSA_cen_sto_McCormick4}
of the temporally aggregated model \eqref{eq:TSA_cen_sto_MPC} yields
\begin{equation}
\sum_{k \in \boldsymbol{K}_r} p^{\mathrm{h}}_{n,\omega,t+k|t}
\geq
10^{-6} \, w \, g \, \eta_n \, \underline{H}_{n} \sum_{k \in \boldsymbol{K}_r} q^{\mathrm{tr}}_{n,\omega,t+k|t},
\; \forall n, \forall \omega, \forall r, 
\label{prop:TSA_cen_sto_McCormick1}
\end{equation}
\begin{multline}
\sum_{k \in \boldsymbol{K}_r} p^{\mathrm{h}}_{n,\omega,t+k|t}
\geq
10^{-6} \, w \, g \, \eta_n \left(
\overline{Q}^{\mathrm{tr}}_{n}
\left(h_{n,\omega,t+\min_{k \in \boldsymbol{K}_r}(t+k)|t} + \left(K_r - 1\right) \underline{H}_n\right)
+ \overline{H}_{n} \sum_{k \in \boldsymbol{K}_r} q^{\mathrm{tr}}_{n, \omega, t+k|t} - K_r \, \overline{Q}^{\mathrm{tr}}_{n}  \, \overline{H}_{n}\right),\\
\forall n, \forall \omega, \forall r, 
\label{prop:TSA_cen_sto_McCormick2}
\end{multline}
\begin{equation}
\sum_{k \in \boldsymbol{K}_r} p^{\mathrm{h}}_{n,\omega,t+k|t}
\leq
10^{-6} \, w \, g \, \eta_n \, \overline{H}_{n} \sum_{k \in \boldsymbol{K}_r} q^{\mathrm{tr}}_{n,\omega,t+k|t},
\; \forall n, \forall \omega, \forall r,
\label{prop:TSA_cen_sto_McCormick3}
\end{equation}
\begin{multline}
\sum_{k \in \boldsymbol{K}_r} p^{\mathrm{h}}_{n,\omega,t+k|t}
\leq
10^{-6} \, w \, g \, \eta_n \left(\overline{Q}^{\mathrm{tr}}_{n}
\left(
h_{n, \omega, \min_{k \in \boldsymbol{K}_r}(t+k)|t} + \left(K_r - 1\right) \overline{H}_n
\right)
+ \underline{H}_{n} \, \sum_{k \in \boldsymbol{K}_r} q^{\mathrm{tr}}_{n, \omega, t+k|t} - K_r \,\overline{Q}^{\mathrm{tr}}_n \, \underline{H}_{n}\right), \\
\forall n, \forall \omega, \forall r.
\label{prop:TSA_cen_sto_McCormick4}
\end{multline}

Similarly to the previously analyzed constraints, the aggregated constraints 
\eqref{prop:TSA_cen_sto_McCormick1} and \eqref{prop:TSA_cen_sto_McCormick3}
are directly implied by the full-scale constraints
\eqref{eq:cen_sto_McCormick1} and \eqref{eq:cen_sto_McCormick3},
respectively, by summing the corresponding full-scale constraints over all time periods $k$ assigned to each cluster $\boldsymbol{K}_r$.
Moreover, from \eqref{eq:cen_sto_McCormick2} it follows that
\begin{equation}
\sum_{k \in \boldsymbol{K}_r} p^{\mathrm{h}}_{n, \omega, t+k|t}
\geq
10^{-6} \, w \, g \, \eta_n
\left(
\overline{Q}^{\mathrm{tr}}_{n} \sum_{k \in \boldsymbol{K}_r} h_{n, \omega, t+k|t} +
\overline{H}_{n} \sum_{k \in \boldsymbol{K}_r} q^{\mathrm{tr}}_{n, \omega, t+k|t} -
K_r \, \overline{Q}^{\mathrm{tr}}_{n} \, \overline{H}_{n}\right),
\; \forall n, \forall \omega, \forall r,
\end{equation}
and
\begin{align*}
& 10^{-6} \, w \, g \, \eta_n
\left(
\overline{Q}^{\mathrm{tr}}_{n} \sum_{k \in \boldsymbol{K}_r} h_{n, \omega, t+k|t} +
\overline{H}_{n} \sum_{k \in \boldsymbol{K}_r} q^{\mathrm{tr}}_{n, \omega, t+k|t} -
K_r \, \overline{Q}^{\mathrm{tr}}_{n} \, \overline{H}_{n}
\right)
= \\
& 10^{-6} \, w \, g \, \eta_n
\left(
\overline{Q}^{\mathrm{tr}}_{n} \left(h_{n, \omega, \min_{k \in \boldsymbol{K}_r}(t+k)|t} + \sum_{k \in \boldsymbol{K}_r \setminus \left\{\min_{k \in \boldsymbol{K}_r}(t+k)\right\}} h_{n, \omega, t+k|t}\right) +
\overline{H}_{n} \sum_{k \in \boldsymbol{K}_r} q^{\mathrm{tr}}_{n, \omega, t+k|t} -
K_r \, \overline{Q}^{\mathrm{tr}}_{n} \, \overline{H}_{n}
\right)
\geq\\
& 10^{-6} \, w \, g \, \eta_n
\left(
\overline{Q}^{\mathrm{tr}}_{n} \left(h_{n, \omega, \min_{k \in \boldsymbol{K}_r}(t+k)|t} + \left(K_r - 1\right) \underline{H}_n\right) +
\overline{H}_{n} \sum_{k \in \boldsymbol{K}_r} q^{\mathrm{tr}}_{n, \omega, t+k|t} -
K_r \, \overline{Q}^{\mathrm{tr}}_{n} \, \overline{H}_{n}
\right),
\; \forall n, \forall \omega, \forall r,
\end{align*}
where the right-hand-side term in the last inequality coincides with the right-hand side term of \eqref{prop:TSA_cen_sto_McCormick2}.
Consequently, every feasible solution of the full-scale model \eqref{eq:cen_sto_MPC} satisfies the aggregated constraints \eqref{prop:TSA_cen_sto_McCormick2}.
Analogously, the full-scale constraints \eqref{eq:cen_sto_McCormick4} imply the aggregated constraints \eqref{prop:TSA_cen_sto_McCormick4}.

Furthermore, from \eqref{prop_agg_vars5}, \eqref{prop_agg_vars6} and \eqref{prop_agg_vars11},
substituting $\boldsymbol{\bar{z}}$ into \eqref{eq:TSA_cen_sto_water_level}
of the temporally aggregated model \eqref{eq:TSA_cen_sto_MPC} yields
\begin{equation}
l_{n,\omega, \min_{k \in \boldsymbol{K}_r}(t+k)|t}
=
l_{n,\omega, \min_{k \in \boldsymbol{K}_{r-1}}(t+k)|t} +
\sum_{k \in \boldsymbol{K}_r}
\left(q^{\mathrm{in}}_{n,\omega,t+k|t} - q^{\mathrm{out}}_{n,\omega,t+k|t}\right)
\frac{\Delta}{S_n}, \; \forall n, \forall \omega, \forall r \in \boldsymbol{R} \setminus \{0\}.
\label{prop:TSA_cen_sto_water_level}
\end{equation}
Under Assumption \textit{A1}, the time periods $k$ in each cluster are consecutive and preserve the original temporal ordering.
Therefore, summing the full-scale constraints \eqref{eq:cen_sto_water_level} over all time periods $k$ belonging to each cluster $\boldsymbol{K}_r$ yields exactly \eqref{prop:TSA_cen_sto_water_level}.
Consequently, every feasible solution of the full-scale model \eqref{eq:cen_sto_MPC} satisfies the aggregated constraints \eqref{prop:TSA_cen_sto_water_level}.

Moreover, from \eqref{prop_agg_vars9} and \eqref{prop_agg_vars11}, substituting $\boldsymbol{\bar{z}}$ into
\eqref{eq:TSA_cen_sto_wl_lim} and \eqref{eq:TSA_cen_sto_head} of the temporally aggregated model
\eqref{eq:TSA_cen_sto_MPC} yields
\begin{align}
\underline{L}_{n} & \leq
l_{n,\omega, \min_{k \in \boldsymbol{K}_r}(t+k)|t}
\leq \overline{L}_{n},
\; \forall n, \forall \omega, \forall r,
\label{prop:TSA_cen_sto_wl_lim}\\
h_{n,\omega, \min_{k \in \boldsymbol{K}_r}(t+k)|t} & =
l_{n,\omega, \min_{k \in \boldsymbol{K}_r}(t+k)|t} - L^{\mathrm{tlr}}_{n},
\; \forall n, \forall \omega, \forall r,
\label{prop:TSA_cen_sto_head}
\end{align}
which are directly implied by the corresponding full-scale constraints \eqref{eq:cen_sto_wl_lim} and \eqref{eq:cen_sto_head}.

From \eqref{prop_agg_vars7}, by substituting $\boldsymbol{\bar{z}}$ into \eqref{eq:TSA_cen_sto_ramping1} of the temporally aggregated model \eqref{eq:TSA_cen_sto_MPC} yields
\begin{equation}
\frac{1}{K_r} \sum_{k \in \boldsymbol{K}_r} q^{\mathrm{tr}}_{n, \omega, t+k|t}
- \sum_{k \in \boldsymbol{K}_{r-1}} q^{\mathrm{tr}}_{n, \omega, t+k|t}
\leq
\Delta^{\mathrm{tr}}_{n} + \Delta^{\mathrm{tr}}_{n} \; \frac{K_r - 1}{2}, \; \forall n, \forall \omega, \forall r \in \boldsymbol{R} \setminus \{0\}.
\label{prop:TSA_cen_sto_ramping1}
\end{equation}
Moreover, from the ramping constraints \eqref{eq:cen_sto_ramping1} and \eqref{eq:cen_sto_ramping2}, the left-hand side of \eqref{prop:TSA_cen_sto_ramping1} can be upper bounded as follows:
\begin{multline}
\frac{1}{K_r} \sum_{k \in \boldsymbol{K}_r} q^{\mathrm{tr}}_{n, \omega, t+k|t}
- \sum_{k \in \boldsymbol{K}_{r-1}} q^{\mathrm{tr}}_{n, \omega, t+k|t}
\leq
q^{\mathrm{tr}}_{n, \omega, \min_{k \in \boldsymbol{K}_r}(t+k)|t} + \Delta^{\mathrm{tr}}_{n} \, \frac{K_r - 1}{2} - q^{\mathrm{tr}}_{n, \omega, \max_{k \in \boldsymbol{K}_{r-1}}(t+k)|t},\\
\forall n, \forall \omega, \forall r \in \boldsymbol{R} \setminus \{0\}.
\label{prop:TSA_cen_sto_ramping1_proof1}
\end{multline}
Under Assumption \textit{A1}, $\min_{k \in \boldsymbol{K}_r}(t+k)$ and $\max_{k \in \boldsymbol{K}_{r-1}}(t+k)$ are consecutive time steps in the full-scale prediction horizon $\boldsymbol{K}$. Therefore, from the full-scale ramping constraints \eqref{eq:cen_sto_ramping1} and \eqref{eq:cen_sto_ramping2}, it follows that:
\begin{equation}
q^{\mathrm{tr}}_{n, \omega, \min_{k \in \boldsymbol{K}_r}(t+k)|t} - q^{\mathrm{tr}}_{n, \omega, \max_{k \in \boldsymbol{K}_{r-1}}(t+k)|t} \leq \Delta^{\mathrm{tr}}_{n}, \; \forall n, \forall \omega, \forall r \in \boldsymbol{R} \setminus \{0\}.
\label{prop:TSA_cen_sto_ramping1_proof2}
\end{equation}
Combining \eqref{prop:TSA_cen_sto_ramping1_proof1} and \eqref{prop:TSA_cen_sto_ramping1_proof2}, which are both implied by the full-scale ramping constraints \eqref{eq:cen_sto_ramping1} and \eqref{eq:cen_sto_ramping2}, it directly follows that every feasible solution of the full-scale model \eqref{eq:cen_sto_MPC} satisfies the aggregated constraints \eqref{prop:TSA_cen_sto_ramping1}.
Analogously, the full-scale constraints \eqref{eq:cen_sto_ramping1} and \eqref{eq:cen_sto_ramping2} imply the aggregated constraints \eqref{eq:TSA_cen_sto_ramping2}.

From \eqref{prop_agg_vars10}, substituting $\boldsymbol{\bar{z}}$ into
\eqref{eq:TSA_cen_sto_binary_lim} of the temporally aggregated model
\eqref{eq:TSA_cen_sto_MPC} yields
\begin{equation}
\frac{1}{K_r} \sum_{k \in \boldsymbol{K}_r}
b_{n,\omega,t+k|t} \in \left\{ \frac{m}{K_r} \,\middle|\, m = 0,1,\dots,K_r \right\}, \; \forall n, \forall \omega, \forall r,
\end{equation}
which directly follows from the full-scale constraints \eqref{eq:cen_sto_binary_lim}.

Finally, the aggregated constraints
\eqref{eq:TSA_cen_sto_water_level_init},
\eqref{eq:TSA_cen_sto_water_level_final},
\eqref{eq:TSA_cen_sto_MPC_control_hydro}, and
\eqref{eq:TSA_cen_sto_MPC_control_vres}
are preserved at full-scale resolution under Assumption \textit{A2}.
Hence, they are equivalent to their respective full-scale counterparts
\eqref{eq:cen_sto_water_level_init},
\eqref{eq:cen_sto_water_level_final},
\eqref{eq:cen_sto_MPC_control_hydro}, and
\eqref{eq:cen_sto_MPC_control_vres}.

It follows that any $\boldsymbol{\bar{z}}$ constructed via \eqref{prop_agg_vars1}--\eqref{prop_agg_vars11} from a feasible solution $\boldsymbol{z}$ of the full-scale centralized stochastic model \eqref{eq:cen_sto_MPC}
constitutes a feasible solution for the temporally aggregated centralized stochastic dispatch model \eqref{eq:TSA_cen_sto_MPC}.

We next demonstrate that, for any feasible solution $\boldsymbol{\bar{z}}$ of the temporally aggregated centralized stochastic dispatch model constructed via
\eqref{prop_agg_vars1}--\eqref{prop_agg_vars11}, and its corresponding full-scale solution $\boldsymbol{z}$,
it holds that $\bar{F}(\boldsymbol{\bar{z}}) \leq F(\boldsymbol{z})$.

From \eqref{prop_agg_vars1}, substituting $\boldsymbol{\bar{z}}$ into the objective function \eqref{eq:TSA_cen_sto_obj} of the temporally aggregated model \eqref{eq:TSA_cen_sto_MPC} yields
\begin{equation}
\bar{F}(\boldsymbol{\bar{z}}) =
\sum_{\omega \in \boldsymbol{\Omega}}
\pi_{\omega}
\sum_{r \in \boldsymbol{R}} K_r
\left(
\frac{1}{K_r} \sum_{k \in \boldsymbol{K}_r}
\left(p_{\omega,t+k|t} - P^{\mathrm{ref}}_{t+k}\right)
\right)^2.
\label{prop:proof_obj1}
\end{equation}

Jensen’s inequality implies that
\begin{equation}
\left(\frac{1}{K_r} \sum_{k \in \boldsymbol{K}_r}
\left(p_{\omega,t+k|t} - P^{\mathrm{ref}}_{t+k}\right)
\right)^2
\leq
\frac{1}{K_r} \sum_{k \in \boldsymbol{K}_r}
\left(p_{\omega,t+k|t} - P^{\mathrm{ref}}_{t+k}\right)^2.
\label{prop:proof_obj2}
\end{equation}

Therefore, combining \eqref{prop:proof_obj1} and \eqref{prop:proof_obj2} yields
\begin{equation}
\bar{F}(\boldsymbol{\bar{z}}) \leq
\sum_{\omega \in \boldsymbol{\Omega}}
\pi_{\omega}
\sum_{r \in \boldsymbol{R}} \sum_{k \in \boldsymbol{K}_r}
\left(p_{\omega,t+k|t} - P^{\mathrm{ref}}_{t+k}\right)^2
=
\sum_{\omega \in \boldsymbol{\Omega}} \pi_{\omega} \sum_{k \in \boldsymbol{K}}
\left(p_{\omega,t+k|t} - P^{\mathrm{ref}}_{t+k}\right)^2
=
F(\boldsymbol{z}),
\end{equation}
where $F(\boldsymbol{z})$ is defined in \eqref{eq:cen_sto_obj}.
\end{proof}

In summary, Proposition~\ref{prop:main_result} establishes that every feasible solution $\boldsymbol{z}$ of the full-scale model \eqref{eq:cen_sto_MPC} can be mapped to a feasible solution $\boldsymbol{\bar{z}}$ of the temporally aggregated model \eqref{eq:TSA_cen_sto_MPC}, with an objective function value that is lower than or equal to that of the original full-scale model.
Consequently, solving the temporally aggregated centralized stochastic MPC scheme introduced in this subsection always yields a lower bound on the optimal objective function value of the full-scale controller presented in Subsection~\ref{subsec:SMPC} at each MPC iteration. 
Importantly, this result is independent of the specific clustering technique used to construct the temporally aggregated model \eqref{eq:TSA_cen_sto_MPC} provided that Assumptions~\textit{A1} and~\textit{A2} of Proposition~\ref{prop:main_result} hold.

\subsection{Temporally Aggregated Distributed Stochastic Model Predictive Control Scheme}
\label{subsec:DSMPC_TSA}

While the temporally aggregated centralized stochastic MPC scheme presented in Subsection~\ref{subsec:SMPC_TSA} reduces the temporal dimensionality of the original full-scale controller introduced in Subsection~\ref{subsec:SMPC}, 
thereby offering substantial potential computational savings as the prediction horizon $\boldsymbol{K}$ increases,
it does not include a dedicated mechanism to address scalability with respect to the number of considered scenarios,
which constitutes another major source of computational complexity.
To address this limitation, this subsection employs ADMM to decompose the temporally aggregated centralized stochastic MPC scheme of Subsection~\ref{subsec:SMPC_TSA} across scenarios.

Notably, the optimization problem \eqref{eq:TSA_cen_sto_MPC} can be naturally interpreted as a \textit{consensus problem},
in which the non-anticipativity constraints \eqref{eq:TSA_cen_sto_MPC_control_hydro} and \eqref{eq:TSA_cen_sto_MPC_control_vres}
act as consensus constraints enforcing consistency of the control actions at the current time period ($k = 0$) across all scenarios.
By exploiting this inherent structural property, consensus ADMM \citep{boyd2010distributed} can be directly applied to decompose \eqref{eq:TSA_cen_sto_MPC} across scenarios.

Let
\begin{equation}
\begin{aligned}
\boldsymbol{\bar{z}}_{\omega} \coloneqq \Big\{
& \bar{p}_{\omega,t+r|t},
\bar{p}^{\mathrm{h}}_{n,\omega,t+r|t},
\bar{p}^{\mathrm{w}}_{\omega,t+r|t},
\bar{p}^{\mathrm{s}}_{\omega,t+r|t},
\bar{l}_{n,\omega,t+r|t},
\bar{q}^{\mathrm{in}}_{n,\omega,t+r|t}, \\
& \bar{q}^{\mathrm{out}}_{n,\omega,t+r|t},
\bar{q}^{\mathrm{tr}}_{n,\omega,t+r|t},
\bar{q}^{\mathrm{br}}_{n,\omega,t+r|t},
\bar{h}_{n,\omega,t+r|t},
\bar{b}_{n,\omega,t+r|t}
\;\Big|\;
n \in \boldsymbol{N},
r \in \boldsymbol{R}
\Big\}, \; \forall \omega,
\end{aligned}
\end{equation}
denote the set of local primal variables associated with scenario $\omega$ in the consensus problem \eqref{eq:TSA_cen_sto_MPC}.
We denote by $\boldsymbol{\Gamma}_\omega(\boldsymbol{\theta}_\omega)$ the feasible set of the variables in $\boldsymbol{\bar{z}}_{\omega}$, defined by the constraints \eqref{eq:TSA_cen_sto_water_level}--\eqref{eq:TSA_cen_sto_power_balance} and \eqref{eq:TSA_cen_sto_binary_lim}, and parameterized by the set of scenario-dependent parameters $\boldsymbol{\theta}_\omega$.
We denote by
$\boldsymbol{\bar{\lambda}}^{\boldsymbol{\mathrm{h}}}_{n,\omega,t|t}$
the vector of dual variables associated with the consensus constraints
\eqref{eq:TSA_cen_sto_MPC_control_hydro} for the hydropower plant $n$ in scenario $\omega$,
and by
$\boldsymbol{\bar{\lambda}}^{\boldsymbol{\mathrm{vRES}}}_{\omega,t|t}$
the vector of dual variables associated with the consensus constraints
\eqref{eq:TSA_cen_sto_MPC_control_vres} in scenario $\omega$.
Furthermore, we denote by $\rho$ the ADMM penalty parameter.

The application of consensus ADMM to \eqref{eq:TSA_cen_sto_MPC} yields the following steps, which are performed over $i \in \boldsymbol{I}$ iterations.

\textbf{Step I}. Local primal variable update:
\begin{align}
\boldsymbol{\bar{z}}^{(i+1)}_\omega \coloneqq \argmin_{\boldsymbol{\bar{z}}_\omega \in \boldsymbol{\Gamma}_\omega(\boldsymbol{\theta}_\omega)}
\Bigg\{
& \pi_{\omega} \sum_{r \in \boldsymbol{R}} K_r \, \left( \bar{p}_{\omega, t+r|t} - \bar{P}^{\mathrm{ref}}_{t+r} \right)^2
+ \sum_{n \in \boldsymbol{N}}
\left(\boldsymbol{\bar{\lambda}}^{\boldsymbol{\mathrm{h}}^{(i)}}_{n,\omega,t|t}\right)^\top
\begin{bmatrix}
\bar{q}^{\mathrm{br}}_{n,\omega,t|t} \\
\bar{q}^{\mathrm{tr}}_{n,\omega,t|t}
\end{bmatrix}
+ \left(\boldsymbol{\bar{\lambda}}^{\boldsymbol{\mathrm{vRES}}^{(i)}}_{\omega,t|t}\right)^\top
\begin{bmatrix}
\bar{p}^{\mathrm{w}}_{\omega,t|t} \\
\bar{p}^{\mathrm{s}}_{\omega,t|t}
\end{bmatrix}
\nonumber\\
&
+ \frac{\rho^{(i)}}{2}
\sum_{n \in \boldsymbol{N}}
\left\|
\begin{bmatrix}
\bar{q}^{\mathrm{br}}_{n,\omega,t|t} \\[0.4em]
\bar{q}^{\mathrm{tr}}_{n,\omega,t|t}
\end{bmatrix}
- \boldsymbol{\bar{u}}^{\boldsymbol{\mathrm{h}}^{(i)}}_{n,t|t}
\right\|_2^2
+ \frac{\rho^{(i)}}{2}
\left\|
\begin{bmatrix}
\bar{p}^{\mathrm{w}}_{\omega,t|t} \\[0.4em]
\bar{p}^{\mathrm{s}}_{\omega,t|t}
\end{bmatrix}
- \boldsymbol{\bar{u}}^{\boldsymbol{\mathrm{vRES}}^{(i)}}_{t|t}
\right\|_2^2
\Bigg\}, \; \forall \omega,
\label{ADMM_step1}
\end{align}
where $\bar{P}^{\mathrm{ref}}_{t+r} = \frac{1}{K_r} \sum_{k \in \boldsymbol{K}_r} P^{\mathrm{ref}}_{t+k}$, for all $r \in \boldsymbol{R}$.

\textbf{Step II}. Global primal variable update:
\begin{align}
\boldsymbol{\bar{u}}^{\boldsymbol{\mathrm{h}}^{(i+1)}}_{n,t|t}
\coloneqq
& \; \frac{1}{|\boldsymbol{\Omega}|}
\sum_{\omega \in \boldsymbol{\Omega}}
\begin{bmatrix}
\bar{q}^{\mathrm{br}^{(i+1)}}_{n,\omega,t|t} \\[0.4em]
\bar{q}^{\mathrm{tr}^{(i+1)}}_{n,\omega,t|t}
\end{bmatrix},
\; \forall n,
\label{ADMM_step2_1}\\
\boldsymbol{\bar{u}}^{\boldsymbol{\mathrm{vRES}}^{(i+1)}}_{t|t}
\coloneqq
& \; \frac{1}{|\boldsymbol{\Omega}|}
\sum_{\omega \in \boldsymbol{\Omega}}
\begin{bmatrix}
\bar{p}^{\mathrm{w}^{(i+1)}}_{\omega,t|t} \\[0.4em]
\bar{p}^{\mathrm{s}^{(i+1)}}_{\omega,t|t}
\end{bmatrix}.
\label{ADMM_step2_2}
\end{align}

\textbf{Step III}. Dual variable update:
\begin{align}
\boldsymbol{\bar{\lambda}}^{\boldsymbol{\mathrm{h}}^{(i+1)}}_{n,\omega,t|t}
\coloneqq
& \; \boldsymbol{\bar{\lambda}}^{\boldsymbol{\mathrm{h}}^{(i)}}_{n,\omega,t|t}
+ \rho^{(i)}
\left( \,
\begin{bmatrix}
\bar{q}^{{\mathrm{br}}^{(i+1)}}_{n,\omega,t|t} \\[0.4em]
\bar{q}^{{\mathrm{tr}}^{(i+1)}}_{n,\omega,t|t}
\end{bmatrix}
- \boldsymbol{\bar{u}}^{\boldsymbol{\mathrm{h}}^{(i+1)}}_{n,t|t}
\right),
\; \forall n, \forall \omega,
\label{ADMM_step3_1}\\
\boldsymbol{\bar{\lambda}}^{\boldsymbol{\mathrm{vRES}}^{(i+1)}}_{\omega,t|t}
\coloneqq
& \; \boldsymbol{\bar{\lambda}}^{\boldsymbol{\mathrm{vRES}}^{(i)}}_{\omega,t|t}
+ \rho^{(i)}
\left( \,
\begin{bmatrix}
\bar{p}^{{\mathrm{w}}^{(i+1)}}_{\omega,t|t} \\[0.4em]
\bar{p}^{{\mathrm{s}}^{(i+1)}}_{\omega,t|t}
\end{bmatrix}
- \boldsymbol{\bar{u}}^{\boldsymbol{\mathrm{vRES}}^{(i+1)}}_{t|t}
\right),
\; \forall \omega.
\label{ADMM_step3_2}
\end{align}

\textbf{Step IV}. Compute the squared norms of the primal and dual residuals, denoted by $\boldsymbol{r^{\mathrm{p}}}$ and $\boldsymbol{r^{\mathrm{d}}}$, respectively:
\begin{align}
\left\|\boldsymbol{r}^{\boldsymbol{\mathrm{p}}^{(i+1)}}\right\|_2^2 & \coloneqq
\sum_{n \in \boldsymbol{N}} \sum_{\omega \in \boldsymbol{\Omega}}
\left\|
\boldsymbol{\bar{u}}^{\boldsymbol{\mathrm{h}}^{(i+1)}}_{n,t|t} - \begin{bmatrix}
\bar{q}^{{\mathrm{br}}^{(i+1)}}_{n,\omega,t|t}\\[0.4em]
\bar{q}^{{\mathrm{tr}}^{(i+1)}}_{n,\omega,t|t}
\end{bmatrix}
\right\|_2^2
+ \sum_{\omega \in \boldsymbol{\Omega}}
\left\|
\boldsymbol{\bar{u}}^{\boldsymbol{\mathrm{vRES}}^{(i+1)}}_{t|t} - \begin{bmatrix}
\bar{p}^{\mathrm{w}^{(i+1)}}_{\omega,t|t} \\[0.4em]
\bar{p}^{\mathrm{s}^{(i+1)}}_{\omega,t|t}
\end{bmatrix}
\right\|_2^2,
\label{ADMM_step4_1}\\
\left\|\boldsymbol{r}^{\boldsymbol{\mathrm{d}}^{(i+1)}}\right\|_2^2 & \coloneqq
\rho^{{(i)}^2}
\left(
\sum_{n \in \boldsymbol{N}}
\left\|
\boldsymbol{\bar{u}}^{\boldsymbol{\mathrm{h}}^{(i+1)}}_{n,t|t} - \boldsymbol{\bar{u}}^{\boldsymbol{\mathrm{h}}^{(i)}}_{n,t|t}
\right\|_2^2
+ \left\|
\boldsymbol{\bar{u}}^{\boldsymbol{\mathrm{vRES}}^{(i+1)}}_{t|t} - \boldsymbol{\bar{u}}^{\boldsymbol{\mathrm{vRES}}^{(i)}}_{t|t}
\right\|_2^2
\right).
\label{ADMM_step4_2}
\end{align}
The consensus ADMM routine terminates when either the maximum number of iterations is reached or the squared norms of both the primal and dual residuals fall below prescribed convergence thresholds,
denoted by $\epsilon^{\mathrm{p,thr}}$ and $\epsilon^{\mathrm{d,thr}}$, respectively.

Furthermore, the squared norms of the primal and dual residuals, computed according to \eqref{ADMM_step4_1} and \eqref{ADMM_step4_2},
are used to adaptively update the penalty parameter $\rho$. 
The parameter is initialized to $\rho^{0}$ and subsequently adjusted according to the strategy discussed in \cite{boyd2010distributed}, using the user-defined parameters $\tau > 1$ and $\mu > 1$.

\textbf{Step V}. Penalty parameter update:
\begin{equation}
\rho^{(i+1)} \coloneqq
\begin{cases}
\tau \rho^{(i)},
& \text{if } \left\|\boldsymbol{r}^{\boldsymbol{\mathrm{p}}^{(i+1)}}\right\|_2^2 > \mu \left\|\boldsymbol{r}^{\boldsymbol{\mathrm{d}}^{(i+1)}}\right\|_2^2,\\[0.4em]
\rho^{(i)} / \tau,
& \text{if } \left\|\boldsymbol{r}^{\boldsymbol{\mathrm{d}}^{(i+1)}}\right\|_2^2 > \mu \left\|\boldsymbol{r}^{\boldsymbol{\mathrm{p}}^{(i+1)}}\right\|_2^2,\\[0.4em]
\rho^{(i)},
& \text{otherwise}.
\end{cases}
\label{ADMM_step5}
\end{equation}

Notably, the convergence of consensus ADMM is generally guaranteed only when all the resulting subproblems are convex \citep{molzahn2017survey}. 
This condition is not satisfied for \eqref{ADMM_step1} due to the presence of the integrality constraints \eqref{eq:TSA_cen_sto_tr_lim}, \eqref{eq:TSA_cen_sto_ph_lim}, and \eqref{eq:TSA_cen_sto_binary_lim}. 
Therefore, to retain a formal performance guarantee and enable a transparent assessment of the solution accuracy achieved by the proposed controller, 
a lower bound on the optimal objective function value of the temporally aggregated model \eqref{eq:TSA_cen_sto_MPC} is derived using consensus ADMM.

Specifically, by introducing into the ADMM routine the following additional optimization step, a Lagrangian relaxation-based lower bound $F^{\mathrm{ADMM,LB}}$ can be computed on the optimal objective function value of \eqref{eq:TSA_cen_sto_MPC}, as formally demonstrated in \cite{gade2016obtaining}.

\textbf{Step VI}. Deriving a lower bound on the optimal objective function value of \eqref{eq:TSA_cen_sto_MPC}:
\begin{align}
F^{{\mathrm{ADMM,LB}}^{(i+1)}}_{\omega} \coloneqq 
\min_{\boldsymbol{\bar{z}}_\omega \in \boldsymbol{\Gamma}_\omega(\boldsymbol{\theta}_\omega)}
\Bigg\{
& \pi_{\omega} \sum_{r \in \boldsymbol{R}} K_r \, \left( \bar{p}_{\omega, t+r|t} - \bar{P}^{\mathrm{ref}}_{t+r} \right)^2
+ \sum_{n \in \boldsymbol{N}}
\left(\boldsymbol{\bar{\lambda}}^{\boldsymbol{\mathrm{h}}^{(i)}}_{n,\omega,t|t}\right)^\top
\begin{bmatrix}
\bar{q}^{\mathrm{br}}_{n,\omega,t|t} \\
\bar{q}^{\mathrm{tr}}_{n,\omega,t|t}
\end{bmatrix}
\nonumber\\
& + \left(\boldsymbol{\bar{\lambda}}^{\boldsymbol{\mathrm{vRES}}^{(i)}}_{\omega,t|t}\right)^\top
\begin{bmatrix}
\bar{p}^{\mathrm{w}}_{\omega,t|t} \\
\bar{p}^{\mathrm{s}}_{\omega,t|t}
\end{bmatrix}
\Bigg\}, \; \forall \omega,
\label{ADMM_step6_1}
\end{align}
\begin{equation}
F^{{\mathrm{ADMM,LB}}^{(i+1)}} \coloneqq \sum_{\omega \in \boldsymbol{\Omega}} F^{{\mathrm{ADMM,LB}}^{(i+1)}}_{\omega}.
\label{ADMM_step6_2}
\end{equation}

Since solving \eqref{ADMM_step6_1}--\eqref{ADMM_step6_2} yields a lower bound on the optimal objective function value of the temporally aggregated model \eqref{eq:TSA_cen_sto_MPC},
which, as established in Proposition~\ref{prop:main_result}, in turn yields a lower bound on the optimal objective function value of the original dispatch model \eqref{eq:cen_sto_MPC},
the resulting \textbf{temporally aggregated distributed stochastic MPC scheme}, which executes the iterative Steps I–VI at each MPC iteration,
provides a valid lower bound on the optimal objective function value of the original dispatch model \eqref{eq:cen_sto_MPC} at every MPC iteration,
i.e., $F^{{\mathrm{ADMM,LB}}}(\boldsymbol{\bar{z}}) \leq \bar{F}(\boldsymbol{\bar{z}}) \leq F(\boldsymbol{z})$.
Notably, this is achieved by solving $|\boldsymbol{\Omega}|$ temporally aggregated subproblems \eqref{ADMM_step1} in parallel,
instead of the original full-scale centralized stochastic model \eqref{eq:cen_sto_MPC},
thereby ensuring simultaneous scalability in both the temporal and scenario dimensions of the dispatch problem.

\subsection{Temporally Aggregated Distributed Stochastic Model Predictive Control Scheme with a Performance Guarantee}
\label{subsec:obj_fun_bounds}

The previous subsections established that the proposed temporally aggregated distributed stochastic MPC scheme presented in Subsection~\ref{subsec:DSMPC_TSA}
always yields a lower bound on the optimal objective function value of the original full-scale centralized MPC scheme introduced in Subsection~\ref{subsec:SMPC} at each MPC iteration. 
In this subsection, we further show how the proposed control scheme can also be used to compute an upper bound on the optimal objective function value of the original full-scale centralized MPC scheme.
By quantifying the relative difference between the objective function bounds at each MPC iteration, a rigorous performance guarantee for the proposed controller can be established, as detailed in the following.

\begin{algorithm}[t]
\caption{Temporally Aggregated Distributed Stochastic Model Predictive Control Scheme with a Performance Guarantee}\label{alg:control_algorithm}
\begin{algorithmic}[1]
\Require \parbox[t]{\dimexpr\linewidth-\algorithmicindent}{
Dispatch model parameters $\Bigl\{ \boldsymbol{\theta}_{\omega} \;\Big|\; \omega \in \boldsymbol{\Omega} \Bigl\}$,
ADMM initial penalty parameter $\rho^0$,
ADMM primal \\ residual threshold $\epsilon^{\mathrm{p,thr}}$,
ADMM dual residual threshold $\epsilon^{\mathrm{d,thr}}$,
ADMM parameters $\tau$ and $\mu$,\\
ADMM maximum number of iterations $\overline{I}$,
optimality threshold $\epsilon^\mathrm{thr}$,
and algorithm maximum \\ number of iterations $\overline{J}$.
}

\Ensure Optimal control decisions $\left\{\boldsymbol{u}^{\boldsymbol{\mathrm{h}}^{\star}}_{n,t|t} \;\Big|\; n \in \boldsymbol{N}\right\}$
and
$\boldsymbol{u}^{\boldsymbol{\mathrm{vRES}}^{\star}}_{t|t}$, a feasible solution $\boldsymbol{z^\mathrm{f}}$, and the objective function bounds ${{F^\mathrm{UB}}}^\star$ and ${{F^\mathrm{LB}}}^\star$.

\State \textit{Higher layer -- Temporal aggregation}
\State \textit{Initialization}:
$j \gets 0$,
$\epsilon^0 \gets +\infty$;

\While{$\epsilon^{(j)} > \epsilon^{\mathrm{thr}}$ and $j \leq \overline{J}$}

\State \parbox[t]{\dimexpr\linewidth-\algorithmicindent}{Assign the time periods $k \in \boldsymbol{K}$ to
$\left\{\boldsymbol{K}_r^{(j)} \;\Big|\; r \in \boldsymbol{R}^{(j)}\right\}$
using any clustering technique, under Assumptions \textit{A1} and \textit{A2} of Proposition \ref{prop:main_result}, and construct the temporally aggregated model~\eqref{eq:TSA_cen_sto_MPC};}

\State \textit{Lower layer -- Scenario decomposition}

\State \textit{Initialization}:
$i \gets 0$,
$\rho^{(i)} \gets \rho^0$,
$\left\{
\boldsymbol{\bar{u}}^{\boldsymbol{\mathrm{h}}^{(i)}}_{n,t|t} \gets \boldsymbol{0}^2,
\boldsymbol{\bar{u}}^{\boldsymbol{\mathrm{vRES}}^{(i)}}_{t|t} \gets \boldsymbol{0}^2,
\boldsymbol{\bar{\lambda}}^{\boldsymbol{\mathrm{h}}^{(i)}}_{n,\omega,t|t} \gets \boldsymbol{0}^2,
\boldsymbol{\bar{\lambda}}^{\boldsymbol{\mathrm{vRES}}^{(i)}}_{\omega,t|t} \gets \boldsymbol{0}^2
\;\Big|\;
n \in \boldsymbol{N}, \omega \in \boldsymbol{\Omega}
\right\}$;

\While{
$\left\|\boldsymbol{r}^{\boldsymbol{\mathrm{p}}^{(i+1)}}\right\|_2^2 > \epsilon^{\mathrm{p,thr}}$
and
$\left\|\boldsymbol{r}^{\boldsymbol{\mathrm{d}}^{(i+1)}}\right\|_2^2 > \epsilon^{\mathrm{d,thr}}$
and
$i \leq \overline{I}$}

\State Step I:
$\left\{\boldsymbol{\bar{z}}^{(i+1)}_\omega \;\Big|\; \omega \in \boldsymbol{\Omega}\right\}\gets$ \eqref{ADMM_step1};
\Comment{In parallel $\forall \omega$}

\State Step II:
$\left\{\boldsymbol{\bar{u}}^{\boldsymbol{\mathrm{h}}^{(i+1)}}_{n,t|t} \;\Big|\; n \in \boldsymbol{N}\right\} \gets$ \eqref{ADMM_step2_1},
$\boldsymbol{\bar{u}}^{\boldsymbol{\mathrm{vRES}}^{(i+1)}}_{t|t} \gets$ \eqref{ADMM_step2_2};

\State Step III:
$\left\{\boldsymbol{\bar{\lambda}}^{\boldsymbol{\mathrm{h}}^{(i+1)}}_{n,\omega,t|t} \;\Big|\; n \in \boldsymbol{N}, \omega \in \boldsymbol{\Omega}\right\} \gets$ \eqref{ADMM_step3_1},
$\left\{\boldsymbol{\bar{\lambda}}^{\boldsymbol{\mathrm{vRES}}^{(i+1)}}_{\omega,t|t} \;\Big|\; \omega \in \boldsymbol{\Omega}\right\} \gets$ \eqref{ADMM_step3_2};
\Comment{In parallel $\forall \omega$}

\State Step IV:
$\left\|\boldsymbol{r}^{\boldsymbol{\mathrm{p}}^{(i+1)}}\right\|_2^2 \gets$ \eqref{ADMM_step4_1},
$\left\|\boldsymbol{r}^{\boldsymbol{\mathrm{d}}^{(i+1)}}\right\|_2^2 \gets$ \eqref{ADMM_step4_2};

\State Step V:
$\rho^{(i+1)} \gets $ \eqref{ADMM_step5};

\State Step VI:
$F^{{\mathrm{ADMM,LB}}^{(i+1)}} \gets$ \eqref{ADMM_step6_1}--\eqref{ADMM_step6_2};
\Comment{In parallel $\forall \omega$}

\State $i \gets i+1$;

\EndWhile

\vspace{0.1cm}
\State $\left\{\boldsymbol{u}^{\boldsymbol{\mathrm{h}}^{(j+1)}}_{n,t|t} \;\Big|\; n \in \boldsymbol{N}\right\}
\gets
\left\{\boldsymbol{\bar{u}}^{\boldsymbol{\mathrm{h}}^{(i)}}_{n,t|t} \;\Big|\; n \in \boldsymbol{N}\right\}$,
$\boldsymbol{u}^{\boldsymbol{\mathrm{vRES}}^{(j+1)}}_{t|t}
\gets
\boldsymbol{\bar{u}}^{\boldsymbol{\mathrm{vRES}}^{(i)}}_{t|t}$;

\State  $F^{{\mathrm{LB}}^{(j+1)}} \gets F^{{\mathrm{ADMM,LB}}^{(i)}}$

\State $\boldsymbol{z}^{\boldsymbol{{\mathrm{f}}}^{(j+1)}}$ and $F^{{\mathrm{UB}}^{(j+1)}} \gets$ \eqref{eq:cen_sto_MPC} with first-stage decisions fixed to
$\left\{\boldsymbol{u}^{\boldsymbol{\mathrm{h}}^{(j)}}_{n,t|t} \;\Big|\; n \in \boldsymbol{N}\right\}$
and
$\boldsymbol{u}^{\boldsymbol{\mathrm{vRES}}^{(j)}}_{t|t}$;
\Comment{In parallel $\forall \omega$}

\State  $\epsilon^{(j+1)} \gets$ \eqref{eq:opt_gap};

\State $j \gets j+1$;

\EndWhile

\vspace{0.1cm}
\State $\left\{\boldsymbol{u}^{\boldsymbol{\mathrm{h}}^{\star}}_{n,t|t} \;\Big|\; n \in \boldsymbol{N}\right\}
\gets
\left\{\boldsymbol{u}^{\boldsymbol{\mathrm{h}}^{(j)}}_{n,t|t} \;\Big|\; n \in \boldsymbol{N}\right\}$,
$\boldsymbol{u}^{\boldsymbol{\mathrm{vRES}}^{\star}}_{t|t}
\gets
\boldsymbol{u}^{\boldsymbol{\mathrm{vRES}}^{(j)}}_{t|t}$,
$\boldsymbol{z^\mathrm{f}} \gets \boldsymbol{z}^{{\boldsymbol{\mathrm{f}}}^{(j)}}$, $F{{^\mathrm{UB}}^\star} \gets F{{^\mathrm{UB}}^{(j)}}$ and $F{{^\mathrm{LB}}^\star} \gets F{{^\mathrm{LB}}^{(j)}}$;

\end{algorithmic}
\end{algorithm}

Once a lower bound, denoted by $F^{\mathrm{LB}}$, on the optimal objective function value $F^\star$ of the full-scale centralized stochastic dispatch model \eqref{eq:cen_sto_MPC} is obtained,
as established in Subsections~\ref{subsec:SMPC_TSA} and \ref{subsec:DSMPC_TSA}, 
an upper bound $F^{\mathrm{UB}}$ can be computed by solving \eqref{eq:cen_sto_MPC} while fixing the first-stage decision variables
$\boldsymbol{u}^{\mathrm{h}}_{n,t|t}$, for all $n \in \boldsymbol{N}$,
and $\boldsymbol{u}^{\mathrm{vRES}}_{t|t}$
to the values obtained from \eqref{ADMM_step2_1} and \eqref{ADMM_step2_2}, respectively.
Since \eqref{eq:cen_sto_MPC} is a two-stage stochastic optimization model,
fixing the first-stage decisions enables the upper bound computation to be performed in parallel across scenarios, analogously to the lower bound computation.
Notably, because this operation corresponds to projecting the solution of the temporally aggregated model--which may, in general, be infeasible for \eqref{eq:cen_sto_MPC}--onto the feasible region of the original full-scale centralized stochastic model \eqref{eq:cen_sto_MPC},
it not only provides an upper bound but also yields a feasible solution for \eqref{eq:cen_sto_MPC}, denoted by $\boldsymbol{z^\mathrm{f}}$.

The proposed \textbf{temporally aggregated distributed stochastic MPC scheme with a performance guarantee} executes Algorithm~\ref{alg:control_algorithm} at each time period $t \in \boldsymbol{T}$.

Algorithm~\ref{alg:control_algorithm} exhibits a hierarchical structure consisting of a higher and a lower layer, performing $j \in \boldsymbol{J}$ and $i \in \boldsymbol{I}$ iterations, respectively.
We define $\overline{J} \coloneqq |\boldsymbol{J}|$ and $\overline{I} \coloneqq |\boldsymbol{I}|$.
At time $t$, the \textbf{higher layer} receives updated forecasts of the uncertainty sources in the dispatch problem,
namely external inflows from the river and its tributaries, and vRES capacity factors,
together with the current observations of the reservoir water levels $L^0_n$, for all $n \in \boldsymbol{N}$.
A user-defined clustering technique, satisfying Assumptions \textit{A1} and \textit{A2} of Proposition~\ref{prop:main_result}, is then employed to construct the temporally aggregated model~\eqref{eq:TSA_cen_sto_MPC}, which is subsequently passed to the \textbf{lower layer}.
Therein, the temporally aggregated model is decomposed into $|\boldsymbol{\Omega}|$ scenario-wise subproblems,
which are solved in parallel via consensus ADMM, as detailed in Subsection~\ref{subsec:DSMPC_TSA}.
Upon convergence of ADMM, as determined by the residual-based stopping criterion of \cite{boyd2010distributed}, the resulting solution of the temporally aggregated model \eqref{eq:TSA_cen_sto_MPC} is used to compute $F^{\mathrm{UB}}$ and $F^{\mathrm{LB}}$.

Let $\epsilon$ denote the relative difference between the objective function bounds $F^{\mathrm{UB}}$ and $F^{\mathrm{LB}}$, defined as follows:
\begin{equation}
\epsilon \coloneqq 100 \, \frac{F^{\mathrm{UB}} - F^{\mathrm{LB}}}{\max\{F^{\mathrm{UB}}, 1\}}.
\label{eq:opt_gap}
\end{equation}

When $\epsilon$ falls below a prescribed threshold $\epsilon^{\mathrm{thr}}$, the algorithm terminates. 
Otherwise, a new temporally aggregated model is constructed,
for example by selecting an alternative clustering technique, refining the clusters obtained in the previous iteration, or increasing the number of clusters employed.

Upon termination, the algorithm provides the optimal control decisions for the hydropower units $n \in \boldsymbol{N}$,
denoted by $\boldsymbol{u}^{\boldsymbol{\mathrm{h}}^{\star}}_{n,t|t}$,
and for the vRES units,
denoted by $\boldsymbol{u}^{\boldsymbol{\mathrm{vRES}}^{\star}}_{t|t}$,
of the CH-vRES hybrid system.
Moreover, it yields the final objective function bounds $F^{\mathrm{UB}\star}$ and $F^{\mathrm{LB}\star}$, as well as a feasible solution $\boldsymbol{z^\mathrm{f}}$ for the original full-scale centralized stochastic dispatch model \eqref{eq:cen_sto_MPC}.

\begin{figure}
  \centering
  \includegraphics[width=\textwidth]{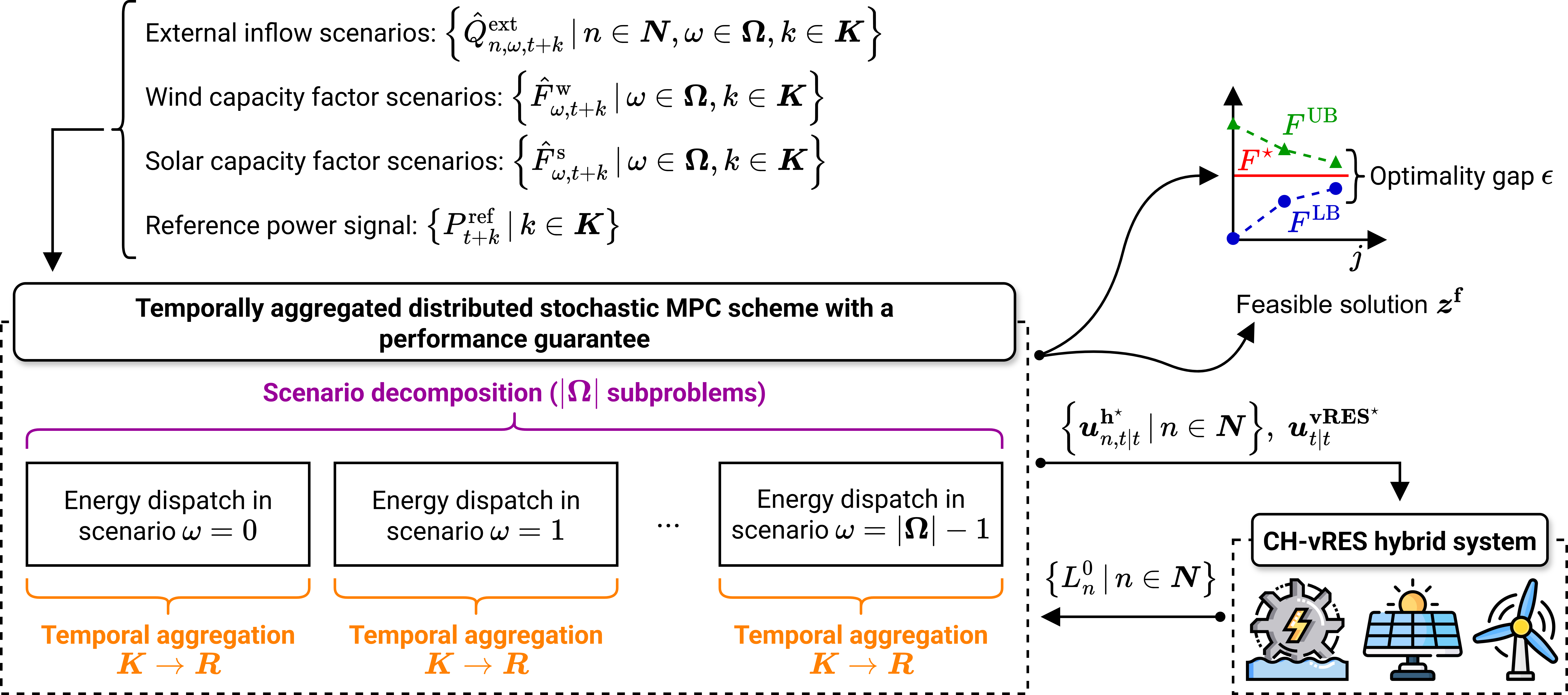}
  \caption{Illustration of the proposed temporally aggregated distributed stochastic MPC scheme with a performance guarantee.}
  \label{fig:fig1_control_algorithm}
\end{figure}

A schematic illustration of Algorithm~\ref{alg:control_algorithm} is provided in Fig.~\ref{fig:fig1_control_algorithm}.
As depicted in the figure, the simultaneous availability of lower and upper bounds on the optimal objective function value of \eqref{eq:cen_sto_MPC} enables the derivation of a formal \textbf{performance guarantee} for the proposed Algorithm~\ref{alg:control_algorithm}.
In particular, the relative difference between the computed objective function bounds provides a rigorous characterization of the achieved \textbf{optimality gap} $\epsilon$,
thereby enabling a theoretically grounded assessment of the algorithm’s accuracy with respect to the original full-scale centralized stochastic dispatch model \eqref{eq:cen_sto_MPC},
while maintaining high computational efficiency through the simultaneous temporal aggregation and scenario-wise parallelization of the dispatch problem.

\section{Simulation Results and Discussion}
\label{sec:results}

This section presents and discusses the simulation results. 
The case study and simulation setup are first introduced in Subsection~\ref{subsec:case_study_description}.
The performance of the proposed control scheme is then evaluated under deterministic and stochastic dispatch formulations in Subsections~\ref{subsec:case_study_deterministic} and~\ref{subsec:case_study_stochastic}, respectively.

\begin{table}
\centering
\renewcommand{\arraystretch}{1}
\begin{tabular}{lll}
\hline
\textbf{Parameter name} & \textbf{Symbol} & \textbf{Value} \\
\hline

Acceleration due to gravity 
& $g$ 
& $9.81\ \mathrm{m/s^{2}}$ \\

Maximum forebay water levels 
& $\overline{L}_{n},\ n \in \{0,1,2\}$ 
& $[123,\ 112,\ 98]\ \mathrm{m}$ \\

Maximum hydropower generation 
& $\overline{P}^{\mathrm{h}}_{n},\ n \in \{0,1,2\}$ 
& $[221,\ 93,\ 136]\ \mathrm{MW}$ \\

Maximum turbine discharges 
& $\overline{Q}^{\mathrm{tr}}_{n},\ n \in \{0,1,2\}$ 
& $[2200,\ 1200,\ 1600]\ \mathrm{m^3/s}$ \\

Minimum barrage discharges 
& $\underline{Q}^{\mathrm{br}}_{n},\ n \in \{0,1,2\}$ 
& $[50,\ 50,\ 50]\ \mathrm{m^3/s}$ \\

Minimum forebay water levels 
& $\underline{L}_{n},\ n \in \{0,1,2\}$ 
& $[120,\ 110,\ 95]\ \mathrm{m}$ \\

Minimum hydropower generation 
& $\underline{P}^{\mathrm{h}}_{n},\ n \in \{0,1,2\}$ 
& $[10,\ 5,\ 11]\ \mathrm{MW}$ \\

Minimum turbine discharges 
& $\underline{Q}^{\mathrm{tr}}_{n},\ n \in \{0,1,2\}$ 
& $[110,\ 60,\ 140]\ \mathrm{m^3/s}$ \\

Ramp limits 
& $\Delta^{\mathrm{tr}}_{n},\ n \in \{0,1,2\}$ 
& $[220,\ 120,\ 160]\ \mathrm{m^3/s}$ \\

Reservoir surface areas 
& $S_{n},\ n \in \{0,1,2\}$ 
& $[1.5,\ 1.2,\ 1.8]\ \mathrm{km^2}$ \\

Solar installed capacity 
& $X^{\mathrm{S}}$ 
& $100\ \mathrm{MW}$ \\

Tailrace water levels 
& $L^{\mathrm{tlr}}_{n},\ n \in \{0,1,2\}$ 
& $[111,\ 101,\ 86]\ \mathrm{m}$ \\

Time period duration 
& $\Delta$ 
& $10\ \mathrm{min}$ \\

Turbine efficiencies 
& $\eta_n,\ n \in \{0,1,2\}$ 
& $[0.9,\ 0.8,\ 0.8]$ \\

Water density 
& $w$ 
& $1000\ \mathrm{kg/m^{3}}$ \\

Wind installed capacity 
& $X^{\mathrm{W}}$ 
& $100\ \mathrm{MW}$ \\
\hline
\end{tabular}
\caption{Parameter values used in the dispatch problem.}
\label{tab:parameters}
\end{table}

\begin{figure}
  \centering
  \includegraphics[width=.71\textwidth]{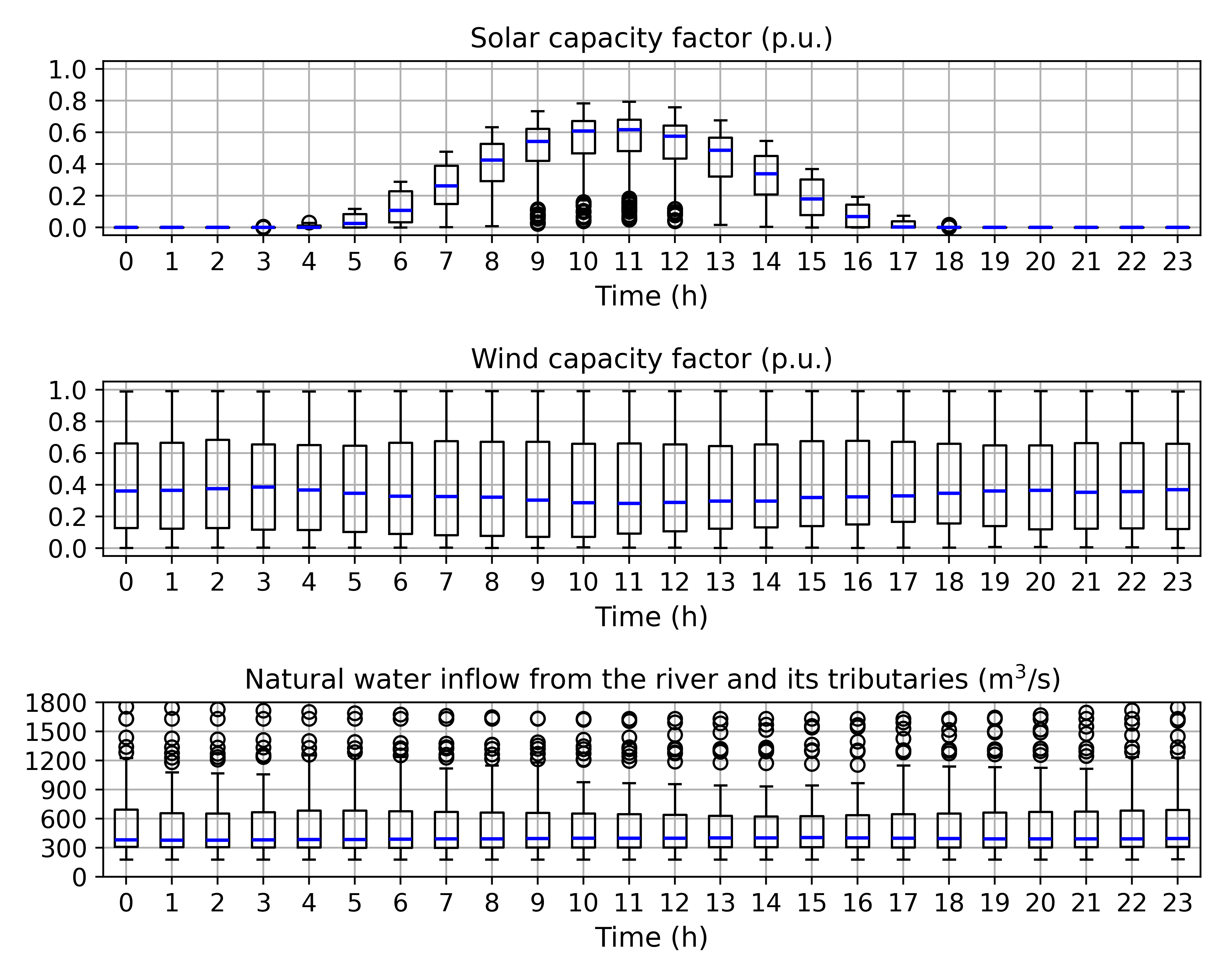}
  \caption{Boxplots of hourly uncertainty realizations. Each boxplot represents the distribution of values observed at a given hour of the day. The box spans the interquartile range, the blue line indicates the median, the whiskers extend to the 10th and 90th percentiles, and outliers are shown as individual points.}
  \label{fig:fig2_observations}
\end{figure}

\begin{figure}
  \centering
  \includegraphics[width=.71\textwidth]{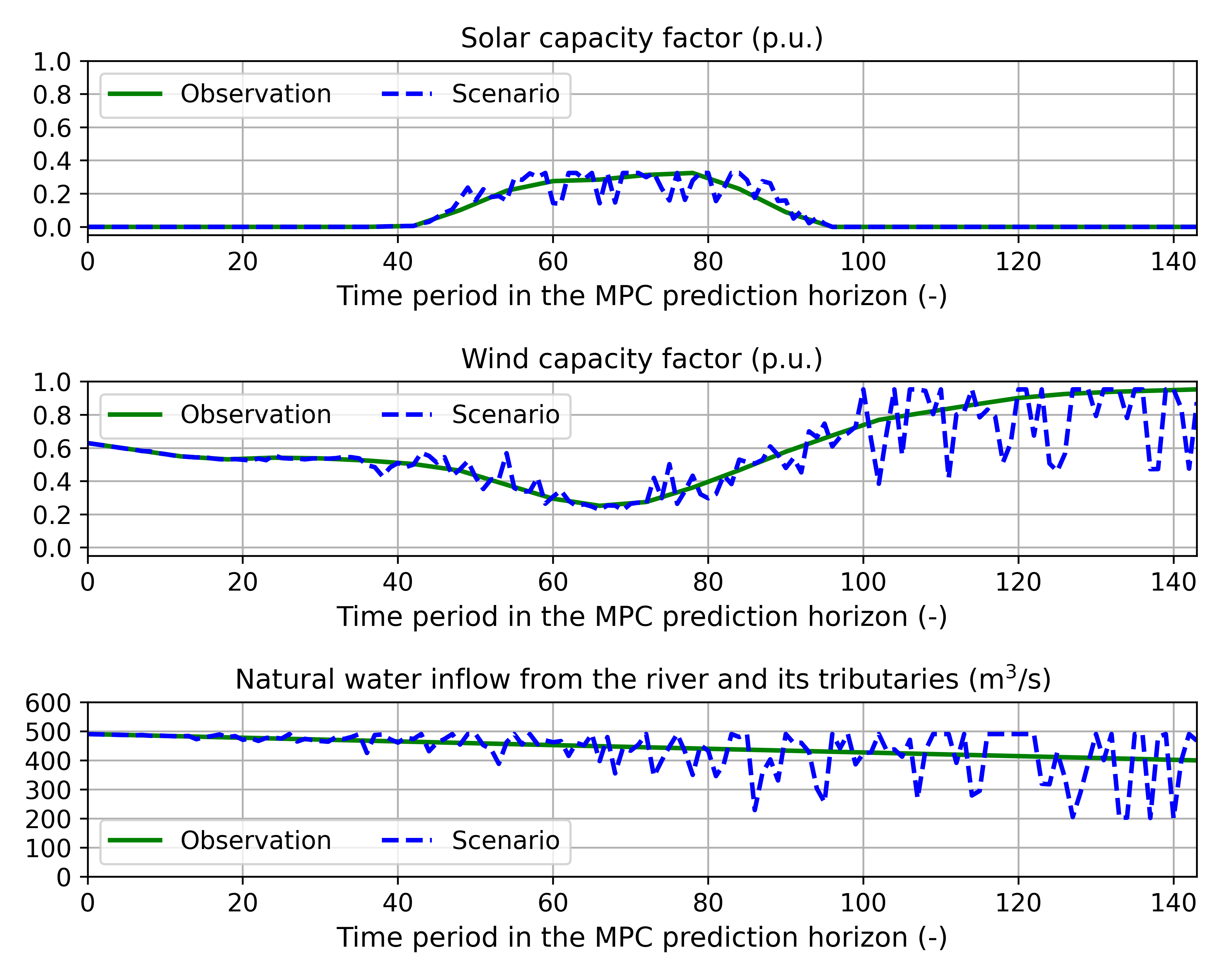}
  \caption{Example of generated uncertainty scenarios at a given iteration of the stochastic MPC scheme, together with the corresponding observed uncertainty realizations.}
  \label{fig:fig3_scenarios_example}
\end{figure}

\subsection{Case Study Description}
\label{subsec:case_study_description}

We consider a case study that mimics a portion of a real-world CH-vRES hybrid system located along the Rhône River in France \citep{piron2016operating}.
The system comprises a cascade of run-of-the-river hydropower plants jointly operated with wind and solar power generation units.
A comprehensive description of the real-world system configuration is provided in \cite{santosuosso2025distributed}.

Specifically, the hydropower cascade is divided into two sections, referred to as the upper Rhône and lower Rhône cascades.
In this study, we focus on the first three hydropower plants of the lower Rhône cascade, as this section is primarily operated to optimize the economic objectives of the hybrid system,
whereas the upper Rhône cascade is mainly used to regulate and stabilize the inflows to the downstream section \citep{santosuosso2025distributed}.
For ease of presentation, the three hydropower plants considered in this study are hereafter denoted, from upstream to downstream, as HPP~0, HPP~1, and HPP~2, respectively.

Unless otherwise specified, the parameters of the CH-vRES hybrid system are set to the values reported in Table~\ref{tab:parameters}, while the parameters of Algorithm~\ref{alg:control_algorithm} are set as follows:
$\rho^0 = 2$,
$\epsilon^{\mathrm{p,thr}} = 10^{-4}$,
$\epsilon^{\mathrm{d,thr}} = 10^{-4}$,
$\tau = 2$, $\mu = 10$,
$\overline{I} = 100$,
$\overline{J} = 100$, and
$\epsilon^{\mathrm{thr}} = 1 \%$.
The forebay water levels of the reservoirs are initialized at their minimum allowable values at the start of the simulation and subsequently updated at each MPC iteration as a function of the computed control actions (see Fig.~\ref{fig:fig1_control_algorithm}).

As discussed in Section~\ref{sec:methodology}, any clustering technique satisfying Assumptions \textit{A1} and \textit{A2} of Proposition~\ref{prop:main_result} can be used in Algorithm~\ref{alg:control_algorithm} without compromising the associated performance guarantee.
We note that standard clustering techniques, such as the widely used k-means algorithm, can be adapted to satisfy Assumptions \textit{A1} and \textit{A2} using, for instance, the approach proposed in \cite{moradi2023capturing}.

In the reported simulations, we adopt the sliding window clustering technique detailed in \cite{santosuosso2026stochastic},
which clusters consecutive time periods $k \in \boldsymbol{K}$ based on whether the distance between their associated clustering features lies within a user-defined similarity threshold.
To progressively tighten the derived objective function bounds, the number of clusters used at each iteration of Algorithm~\ref{alg:control_algorithm} is increased according to the refinement criterion detailed in \cite{santosuosso2026stochastic}.

The proposed controller is evaluated through simulations covering the entire year 2017,
considering a 24-hour prediction horizon and a 10-minute sampling time, yielding $K = 144$.
Several profiles for the power reference signal to be tracked by the hybrid system output are considered, as detailed in the following subsections.

Furthermore, the dispatch problem is subject to three key sources of uncertainty, namely external water inflows into the hydropower plants of the CH-vRES system from the Rhône River and its tributaries, as well as wind and solar capacity factors.
Historical time series representing the observed inflow uncertainty realizations are obtained from the Global Runoff Data Centre \citep{farber2025grdc},
while wind and solar capacity factor data are sourced from Renewables.ninja \citep{pfenninger2016long}.
A characterization of the time series used as the observed uncertainty realizations is provided in Fig.~\ref{fig:fig2_observations}.

Unless otherwise specified, a new set of 20 equiprobable scenarios is generated at each iteration of the stochastic MPC scheme.
The scenarios are obtained by perturbing the time series representing the observed uncertainty realizations with additive Laplace-distributed noise, $\mathcal{L}(0, b_{t+k})$,
while enforcing the physical bounds associated with each uncertainty source
(e.g., wind and solar capacity factors are constrained to the interval $[0,1]$).

Specifically, let $a_{t+k}$ denote the observed realization at time $t+k$.
The scenarios are then generated as $a_{t+k} + \zeta_{t+k}$, where $\zeta_{t+k} \sim \mathcal{L}(0, b_{t+k})$.
The scale parameter $b_{t+k}$ of the Laplace distribution is defined as a function of the observed realization,
with a time-dependent growth factor, as follows:
\begin{equation}
b_{t+k} = \frac{a_{t+k}}{2} \left(\frac{k}{K}\right)^2, \; \forall k \in \boldsymbol{K},
\end{equation}
thereby representing the progressive degradation of forecast accuracy that is typically observed toward the end of the MPC prediction horizon.
An example of the generated scenarios at a given MPC iteration is shown in Fig.~\ref{fig:fig3_scenarios_example}.

\begin{figure}
  \centering
  \includegraphics[width=.72\textwidth]{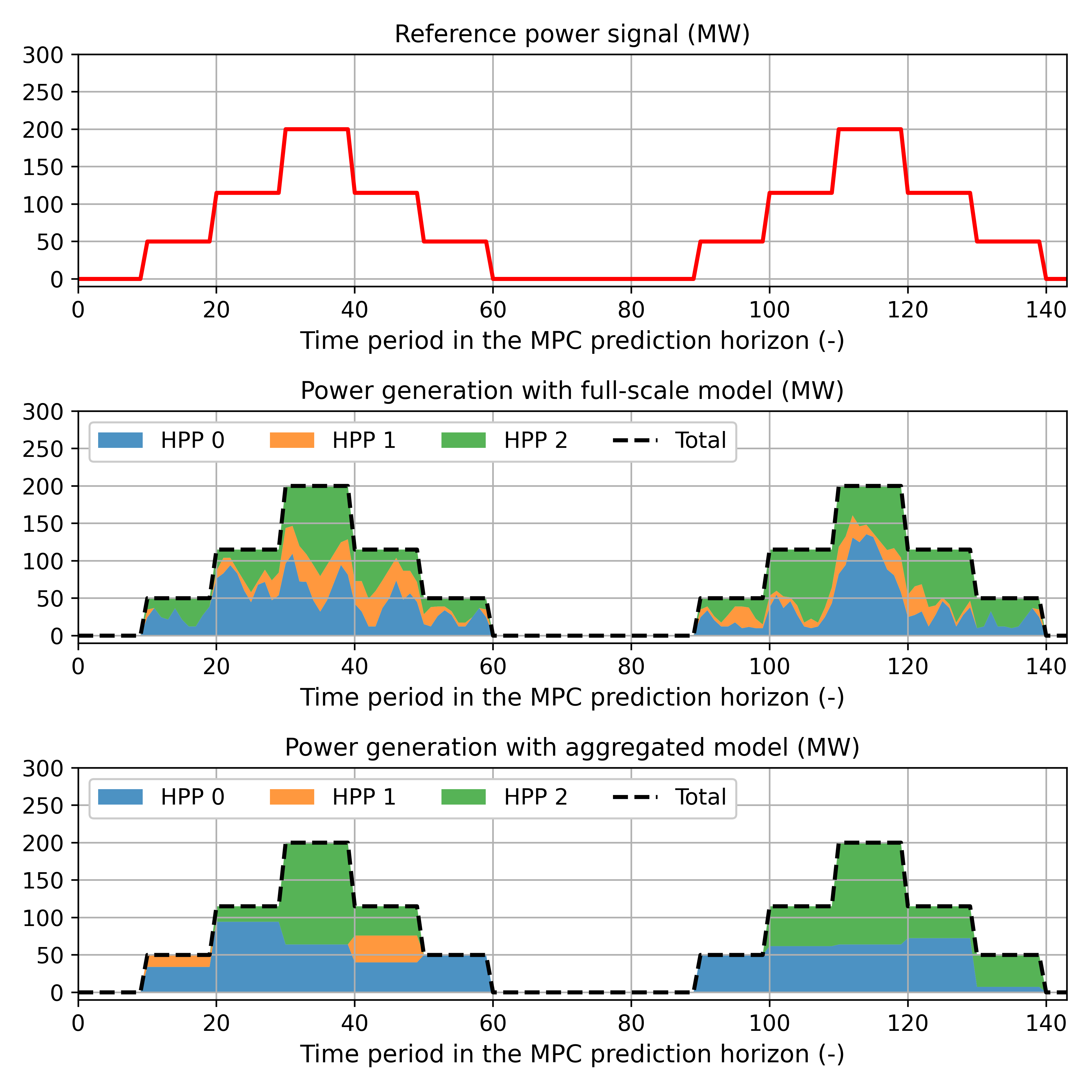}
  \caption{Example of hydropower generation from the cascaded hydropower plants under investigation when using the full-scale dispatch model \eqref{eq:cen_sto_MPC} and the temporally aggregated dispatch model \eqref{eq:TSA_cen_sto_MPC} solved via the proposed Algorithm~\ref{alg:control_algorithm}.}
  \label{fig:fig4_power_output_example1}
\end{figure}

All simulations are conducted on an Intel i7 CPU with 32 GB RAM using Gurobi 12.0.1.
Moreover, Algorithm~\ref{alg:control_algorithm} is implemented using Python's multiprocessing module \cite{multiprocessing}, enabling the derived scenario-wise subproblems to be solved in parallel.

\subsection{Deterministic Dispatch of Cascaded Hydropower Plants}
\label{subsec:case_study_deterministic}

\begin{figure}
  \centering
  \includegraphics[width=.72\textwidth]{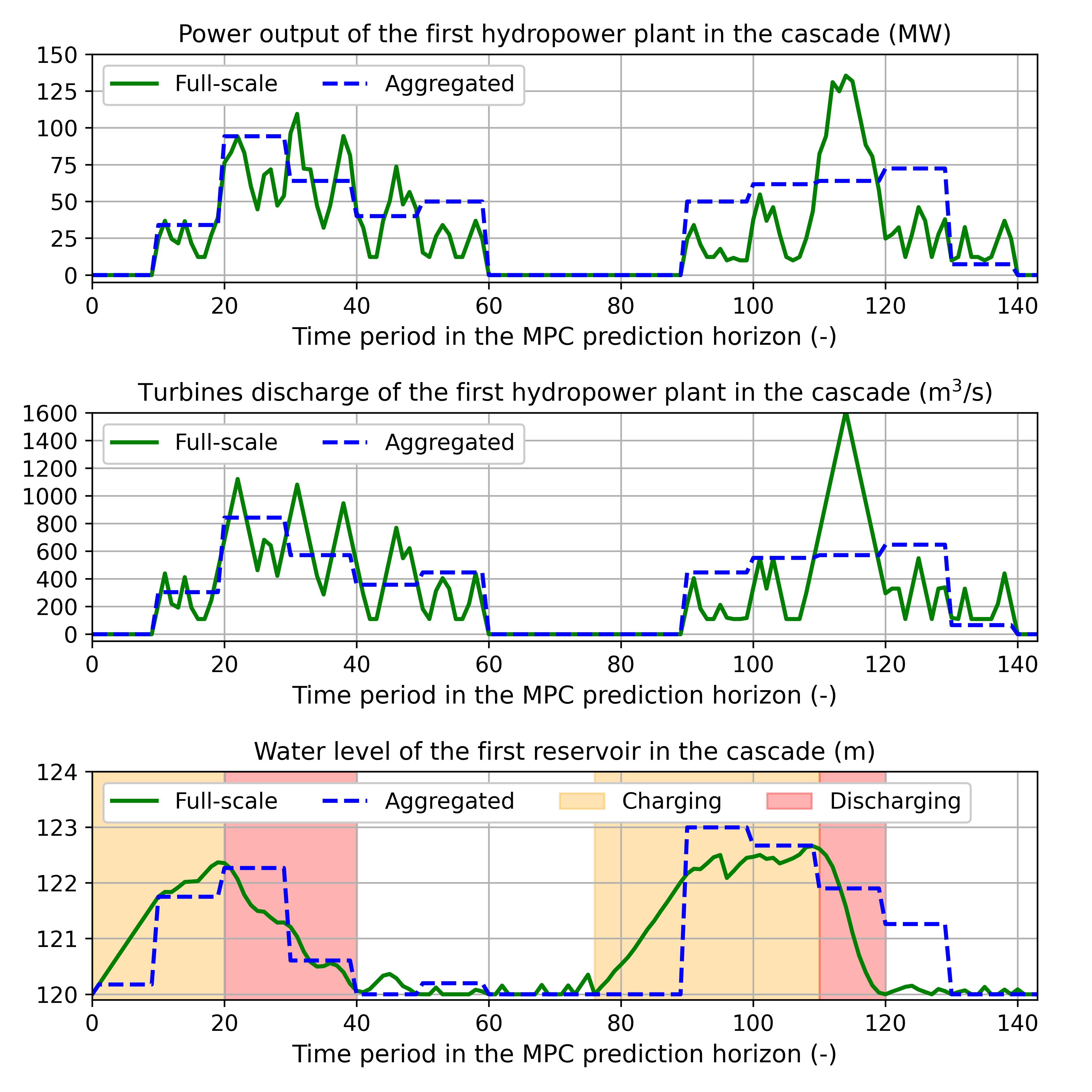}
  \caption{Example of hydropower generation, turbine discharges, and reservoir forebay water level in the first cascaded hydropower plant when employing the full-scale dispatch model \eqref{eq:cen_sto_MPC} and the temporally aggregated dispatch model \eqref{eq:TSA_cen_sto_MPC} solved via the proposed Algorithm~\ref{alg:control_algorithm}.}
  \label{fig:fig5_decisions_example}
\end{figure}

\begin{figure}
  \centering
  \includegraphics[width=.72\textwidth]{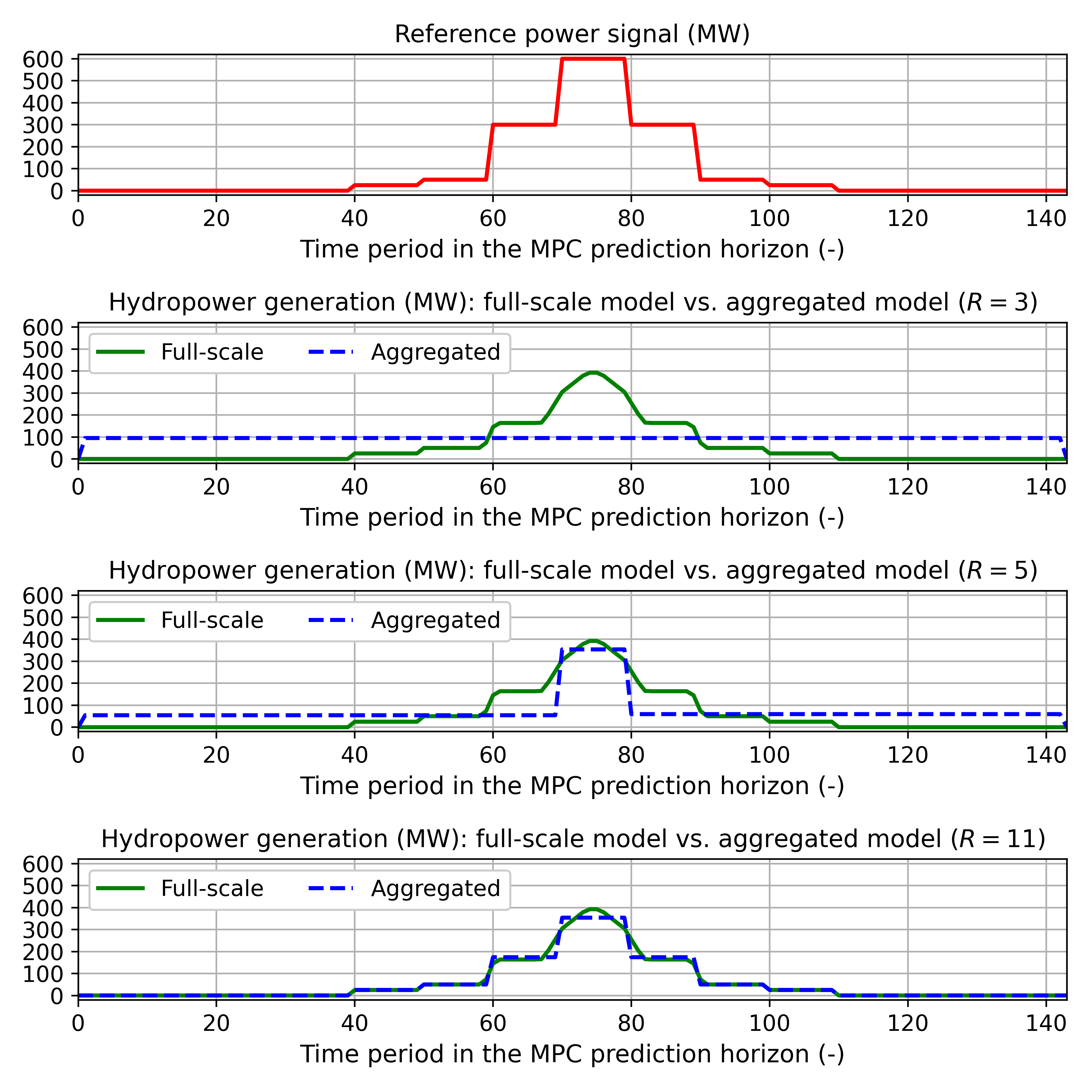}
  \caption{Example of hydropower generation observed when using the temporally aggregated dispatch model \eqref{eq:TSA_cen_sto_MPC} solved via the proposed Algorithm~\ref{alg:control_algorithm}, in comparison with that obtained using the full-scale model \eqref{eq:cen_sto_MPC}, as the number of representative time periods $R$ used for TSA increases.}
  \label{fig:fig6_clusters_sensitivity}
\end{figure}

In this subsection, we consider a deterministic formulation of the energy dispatch problem, corresponding to the MPC problem \eqref{eq:cen_sto_MPC} with $|\boldsymbol{\Omega}| = 1$. 
Since no scenario-wise decomposition is performed in this setting, the simulations reported herein are intended to isolate and assess the impact of the temporal aggregation procedure in the proposed Algorithm~\ref{alg:control_algorithm}. 
Furthermore, to specifically evaluate the effect of temporal aggregation on the primary source of flexibility within the considered hybrid system, namely the cascaded hydropower plants, we set $X^{\mathrm{S}} = 0$ and $X^{\mathrm{W}} = 0$.

Fig.~\ref{fig:fig4_power_output_example1} illustrates an example of power generation obtained using the full-scale MPC scheme of Subsection~\ref{subsec:SMPC},
in comparison with the temporally aggregated MPC scheme of Subsection~\ref{subsec:SMPC_TSA} solved via the proposed Algorithm~\ref{alg:control_algorithm}. Fig.~\ref{fig:fig5_decisions_example} further reports the corresponding power output, turbine discharges, and reservoir forebay water level for the first hydropower plant in the cascade, i.e., HPP~0.
In both figures, the MPC prediction horizon begins at midnight, and the reported results illustrate a typical daily dispatch obtained with the considered MPC schemes.
The reference power signal $P^{\mathrm{ref}}_{t+k}$ is used as the clustering feature within the adopted sliding window clustering technique.
In this case, this yields 15 representative time periods for the temporally aggregated model,
corresponding to the 13 segments of the reference power profile in Fig.~\ref{fig:fig4_power_output_example1}, 
together with two singleton representative time periods associated with the first and last time periods of the MPC prediction horizon,
as required by Assumption~\textit{A2} of Proposition~\ref{prop:main_result}.

As shown in Fig.~\ref{fig:fig4_power_output_example1}, the temporally aggregated model enables the total cascaded hydropower generation to track the reference power signal exactly, as its full-scale counterpart,
while achieving a substantial temporal dimensionality reduction of approximately 90\% (from $K=144$ to $R=15$).
This is obtained despite the fact that the individual generation profiles of the three hydropower plants differ significantly between the full-scale and temporally aggregated models.
This behavior is consistent with both the results shown in Fig.~\ref{fig:fig5_decisions_example} and the theoretical analysis provided in Subsection~\ref{subsec:SMPC_TSA},
according to which the temporally aggregated model does not generally recover an optimal solution of the full-scale model,
but rather produces aggregated (average-like) decisions over the representative time periods.

In particular, as shown in Fig.~\ref{fig:fig5_decisions_example}, both dispatch models exhibit the same qualitative control strategy,
namely the storage of water prior to demand peaks (charging phase) and its subsequent release during demand peaks (discharging phase) to track the reference power signal shown in Fig.~\ref{fig:fig4_power_output_example1}.
However, when the aggregated decisions are disaggregated to the full temporal resolution as shown in Fig.~\ref{fig:fig5_decisions_example}, they may result infeasible for the original full-scale dispatch model.
This observation underscores the role of the projection step introduced in Algorithm~\ref{alg:control_algorithm}, which, as discussed in Subsection~\ref{subsec:obj_fun_bounds},
not only provides an upper bound on the optimal objective function value of the full-scale model,
but also restores feasibility of the decisions derived from the temporally aggregated model.
Notably, when the computed lower bound on the optimal objective function value of the full-scale model is tight, as in the present case,
the feasible solution recovered by Algorithm~\ref{alg:control_algorithm} also typically reproduces an optimal solution of the full-scale model, as further illustrated in the following results.

While the validity of the objective function bounds computed by Algorithm~\ref{alg:control_algorithm} is always guaranteed,
as formally established in Section~\ref{sec:methodology}, 
their tightness, as well as the accuracy of the associated control actions relative to those of the original full-scale model,
depends on the quality of the clustering technique employed for TSA.
For instance, as illustrated in Fig.~\ref{fig:fig6_clusters_sensitivity}, when using the same clustering technique adopted for the previously presented results and the reference power signal as clustering feature, 
the accuracy of the aggregated decisions obtained via Algorithm~\ref{alg:control_algorithm} progressively improves as the number of representative time periods increases (from $R=3$ to $R=11$ in the figure).

\begin{figure}
  \centering
  \includegraphics[width=.73\textwidth]{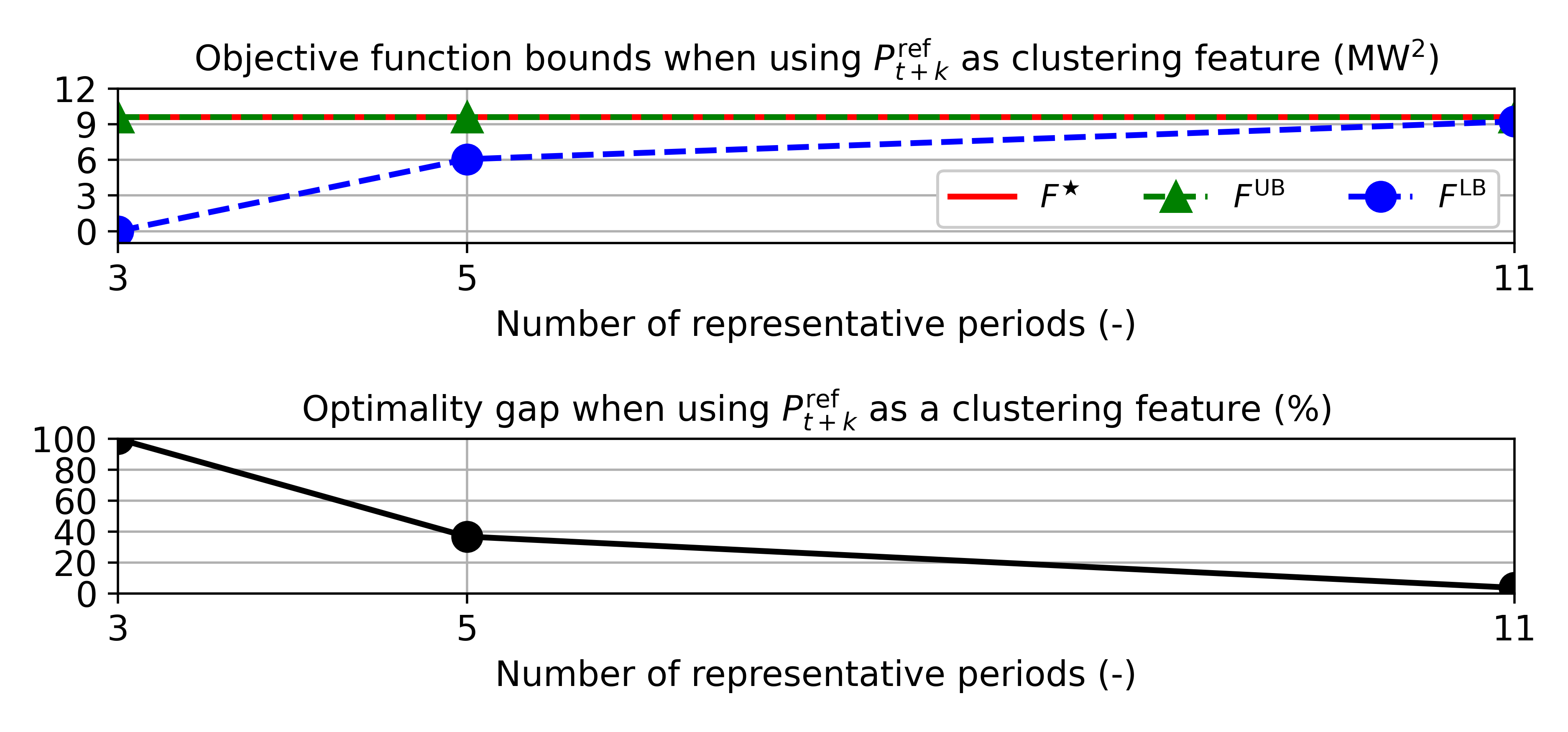}
  \caption{Example of objective function bounds obtained using the proposed Algorithm~\ref{alg:control_algorithm} when the reference power signal is used as clustering feature.}
  \label{fig:fig7_bounds_example1}
\end{figure}

This improvement is clearly reflected in the corresponding objective function bounds computed across the iterations of the proposed control algorithm, as shown in Fig.~\ref{fig:fig7_bounds_example1}.
In particular, an optimality gap below $1\%$, i.e., a solution accuracy exceeding $99\%$, is progressively achieved as the number of representative time periods is increased up to $R=11$.
Notably, as discussed above, the tightness of the computed bounds also provides a direct indication of the quality of the feasible solution obtained alongside the objective function upper bound.
In this case, achieving more than $99\%$ accuracy at $R=11$ informs the decision-maker that
the proposed control algorithm yields a feasible solution for the dispatch problem with an associated objective function error of no more than $1\%$
relative to the optimal objective function value of the original full-scale controller,
while still benefiting from a substantial reduction in computational complexity,
as the dispatch problem is solved over only 11 representative time periods instead of the original 144.

\begin{figure}
  \centering
  \includegraphics[width=.73\textwidth]{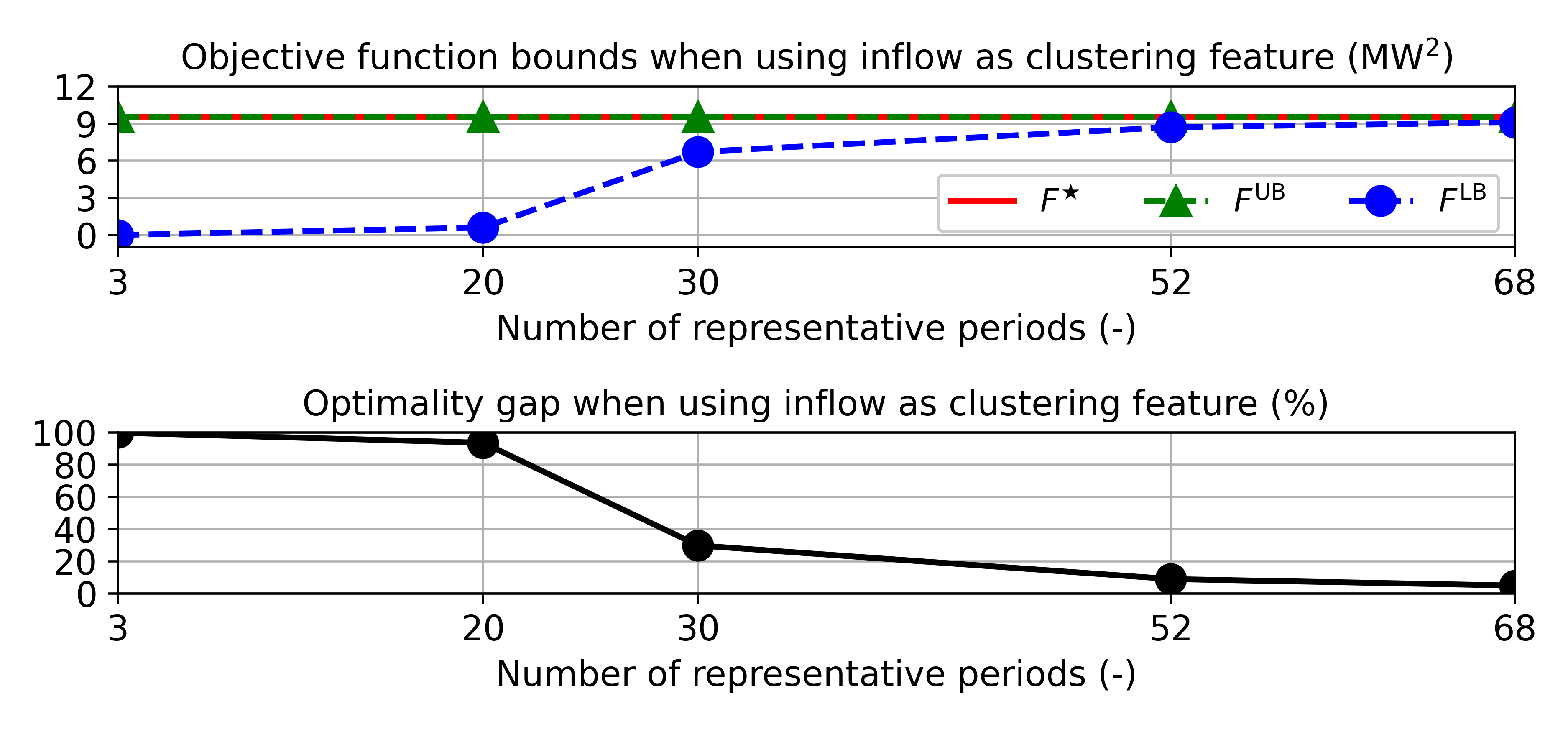}
  \caption{Example of objective function bounds obtained using the proposed Algorithm~\ref{alg:control_algorithm} when the natural water inflow scenario is used as clustering feature.}
  \label{fig:fig8_bounds_example2}
\end{figure}

To further illustrate the sensitivity of Algorithm~\ref{alg:control_algorithm} to the choice of the clustering feature employed in the TSA procedure, 
Fig.~\ref{fig:fig8_bounds_example2} reports the objective function bounds obtained under the same simulation settings as in Figs.~\ref{fig:fig6_clusters_sensitivity} and~\ref{fig:fig7_bounds_example1},
but using the natural water inflow scenario as the clustering feature instead of the reference power signal.
In this case, convergence of the objective function bounds to the desired accuracy is achieved with $R=68$ clusters.

Although this clustering choice requires a larger number of representative time periods compared to the previous case,
the proposed algorithm still achieves a $53\%$ reduction in the temporal dimensionality of the dispatch problem.
In more complex instances of the problem at hand, such a reduction translates into substantial computational savings compared to conventional full-scale MPC schemes, as discussed in the following subsection.

\subsection{Stochastic Dispatch of Cascaded Hydropower Plants and Renewables}
\label{subsec:case_study_stochastic}

In this subsection, we consider the stochastic CH-vRES dispatch problem introduced in Section~\ref{sec:methodology},
in which cascaded hydropower plants, wind, and solar units are jointly operated under the simulation settings described in Subsection~\ref{subsec:case_study_description}.
Unless otherwise specified, the reference power signal $P^{\mathrm{ref}}_{t+k}$ is used as the clustering feature within the adopted sliding window clustering technique.

\begin{figure}
  \centering
  \includegraphics[width=.73\textwidth]{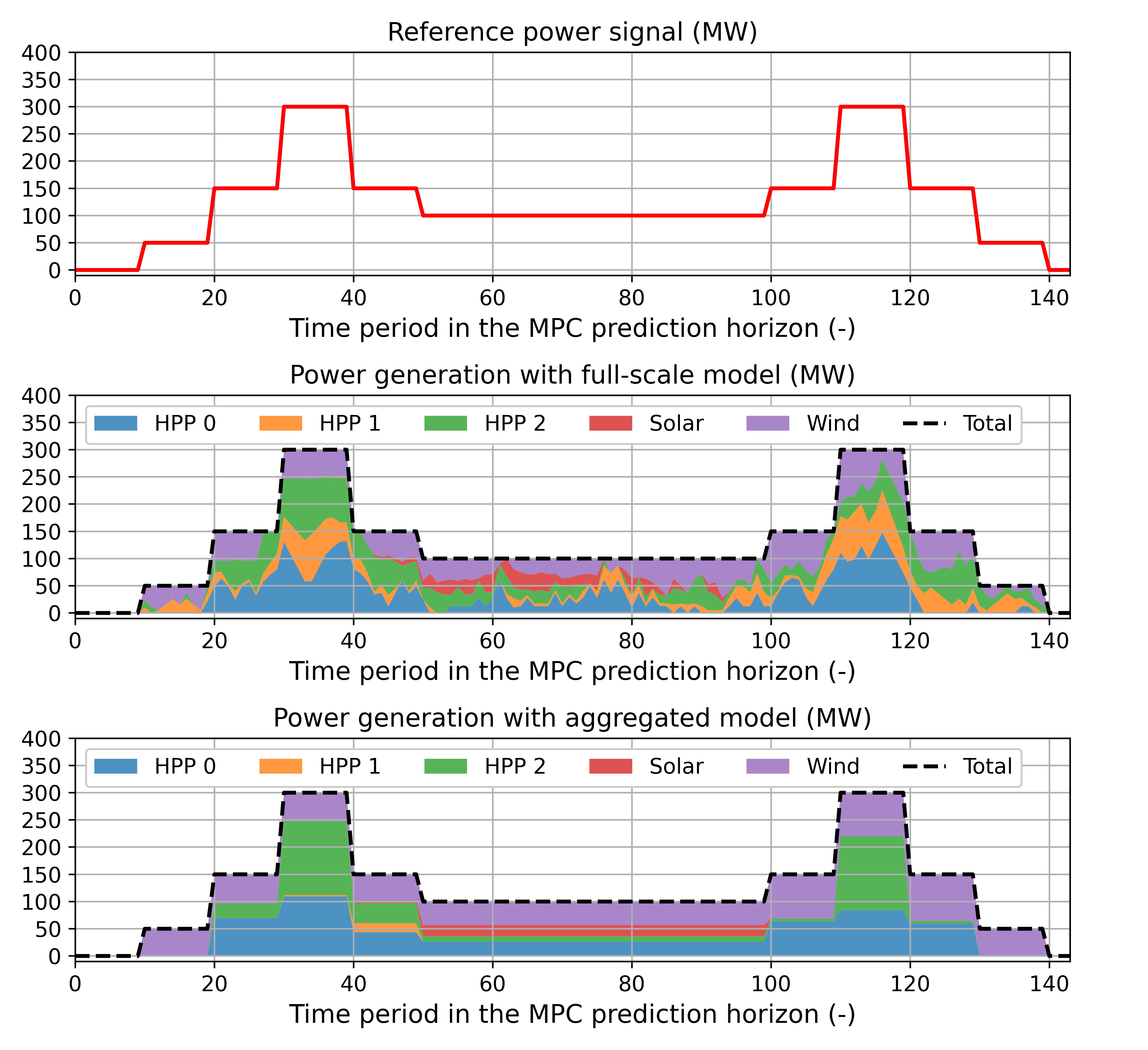}
  \caption{Example of CH-vRES hybrid system power generation when using the full-scale centralized stochastic dispatch model \eqref{eq:cen_sto_MPC} and the temporally aggregated dispatch model \eqref{eq:TSA_cen_sto_MPC} solved via the proposed Algorithm~\ref{alg:control_algorithm}.}
  \label{fig:fig9_power_gen_example}
\end{figure}

Fig.~\ref{fig:fig9_power_gen_example} reports an example of CH-vRES hybrid system dispatch
obtained with the traditional full-scale stochastic centralized MPC scheme of Subsection~\ref{subsec:SMPC},
in comparison with that obtained using the proposed Algorithm~\ref{alg:control_algorithm}.
In this case, the temporally aggregated model employs 13 representative time periods, corresponding to the 11 segments of the reference power profile in Fig.~\ref{fig:fig4_power_output_example1},
together with two singleton representative time periods associated with the first and last time steps of the MPC prediction horizon,
as required by Assumption~\textit{A2} of Proposition~\ref{prop:main_result}.
The MPC prediction horizon begins at midnight, and the reported results illustrate a typical daily dispatch obtained with the considered MPC schemes.

As shown in Fig.~\ref{fig:fig9_power_gen_example}, the temporally aggregated model enables exact tracking of the reference power signal,
thereby recovering the same optimal objective function value of its full-scale counterpart.
Consistently with the deterministic case illustrated in Fig.~\ref{fig:fig4_power_output_example1},
the individual control actions computed for each unit of the CH-vRES hybrid system using the temporally aggregated model are not necessarily feasible for the original full-scale dispatch problem.
Feasibility is therefore restored through the projection step of the proposed Algorithm~\ref{alg:control_algorithm}, as detailed in Subsection~\ref{subsec:obj_fun_bounds}.

\begin{figure}
  \centering
  \includegraphics[width=.73\textwidth]{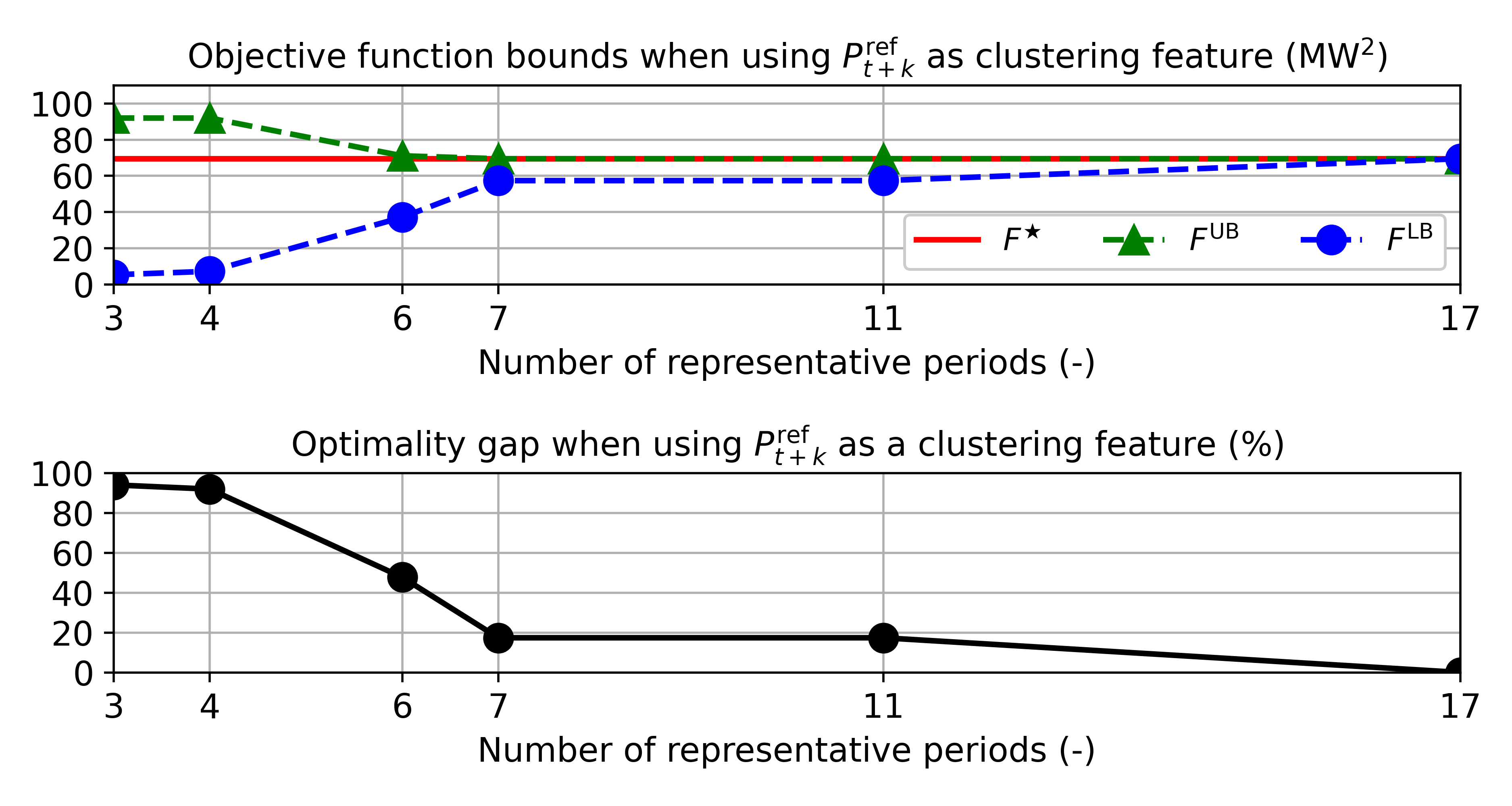}
  \caption{Example of objective function bounds obtained using the proposed Algorithm~\ref{alg:control_algorithm} for the stochastic CH-vRES dispatch problem.}
  \label{fig:fig10_bounds_example_sto}
\end{figure}

An example of the objective function bounds computed by the proposed Algorithm~\ref{alg:control_algorithm} for the stochastic CH-vRES dispatch problem is reported in Fig.~\ref{fig:fig10_bounds_example_sto}.
In contrast to the objective function bounds obtained for the deterministic dispatch problem analyzed in the previous subsection and reported in Figs.~\ref{fig:fig7_bounds_example1} and~\ref{fig:fig8_bounds_example2},
the upper bound does not immediately converge to the optimal objective function value of the full-scale dispatch model. 
This indicates that the proposed control algorithm initially computes a feasible but suboptimal solution for the full-scale dispatch problem.
As the algorithm progressively refines the underlying temporal aggregation, the upper and lower bounds converge, and the associated feasible solution approaches an optimal solution of the full-scale dispatch problem.

\begin{figure}
  \centering
  \includegraphics[width=.73\textwidth]{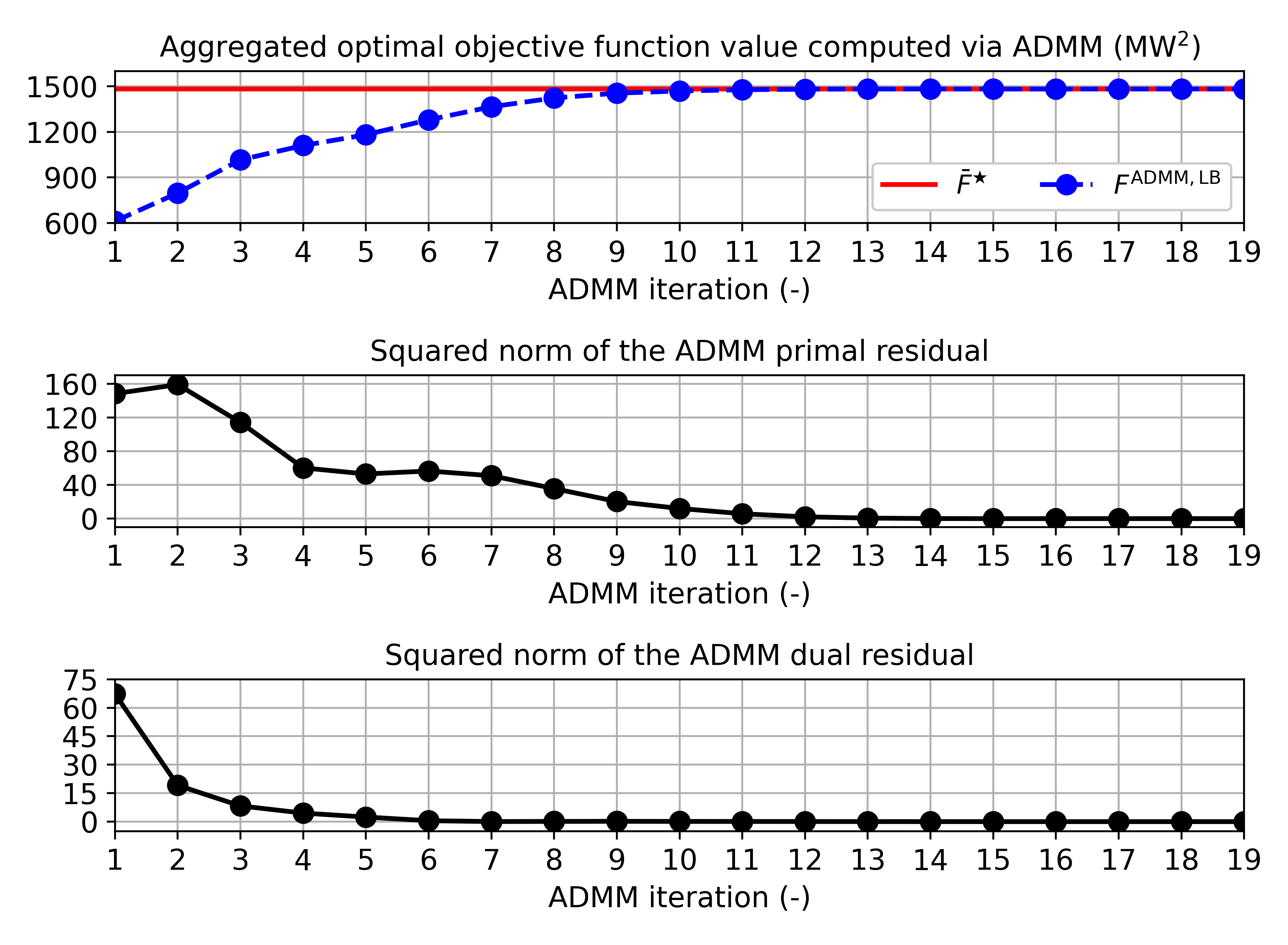}
  \caption{Example of execution of the consensus ADMM routine within the lower layer of the proposed Algorithm~\ref{alg:control_algorithm}.}
  \label{fig:fig11_ADMM_example}
\end{figure}

As detailed in Section~\ref{sec:methodology}, in addition to performing temporal aggregation, Algorithm~\ref{alg:control_algorithm} employs consensus ADMM to decompose the temporally aggregated dispatch model \eqref{eq:TSA_cen_sto_MPC} across scenarios.
An example of the execution of the consensus ADMM routine within Algorithm~\ref{alg:control_algorithm} is reported in Fig.~\ref{fig:fig11_ADMM_example}.
The figure illustrates how the lower bound $F^{\mathrm{ADMM,LB}}$ on the optimal objective function value of the temporally aggregated model, $\bar{F}^{\star}$, is progressively refined throughout the ADMM iterations until convergence is achieved.
In the reported example, convergence is attained after 19 iterations according to the residual-based stopping criterion detailed in \cite{boyd2010distributed}.

\begin{table}
\centering
\renewcommand{\arraystretch}{1}
\begin{tabular}{cccccccc}
\hline
\multirow{2}{*}{\textbf{$|\boldsymbol{K}|$}}
& \multirow{2}{*}{\textbf{$|\boldsymbol{\Omega}|$}}
& \textbf{FS-MPC}
& \textbf{TSA-MPC}
& \textbf{ADMM-MPC}
& \multicolumn{3}{c}{\textbf{Algorithm~\ref{alg:control_algorithm}}}\\
&
& \textbf{Avg. runtime} (s)
& \textbf{Avg. runtime} (s)
& \textbf{Avg. runtime} (s)
& $\epsilon^{\mathrm{thr}}$
& \textbf{Avg.} $R$
& \textbf{Avg. runtime} (s) \\
\hline

\multirow{2}{*}{$144$}
& \multirow{2}{*}{$3$}
& \multirow{2}{*}{$80.2$}
& \multirow{2}{*}{$72.8$ ($\boldsymbol{-9\%}$)}
& \multirow{2}{*}{$81.8$ ($\boldsymbol{+2\%}$)}
& $1\%$
& $91.2$
& $56.7$ ($\boldsymbol{-29\%}$) \\
& & & &
& $5\%$
& $88.1$
& $36.9$ ($\boldsymbol{-54\%}$) \\

\multirow{2}{*}{$288$}
& \multirow{2}{*}{$3$}
& \multirow{2}{*}{$442.4$}
& \multirow{2}{*}{$326.8$ ($\boldsymbol{-26\%}$)}
& \multirow{2}{*}{$430.2$ ($\boldsymbol{-3\%}$)}
& $1\%$
& $101.3$
& $202.6$ ($\boldsymbol{-54\%}$) \\
& & & &
& $5\%$
& $96.8$
& $121.8$ ($\boldsymbol{-72\%}$) \\

\multirow{2}{*}{$144$}
& \multirow{2}{*}{$5$}
& \multirow{2}{*}{$437.8$}
& \multirow{2}{*}{$343.3$ ($\boldsymbol{-22\%}$)}
& \multirow{2}{*}{$408.4$ ($\boldsymbol{-7\%}$)}
& $1\%$
& $93.4$
& $201.3$ ($\boldsymbol{-54\%}$) \\
& & & &
& $5\%$
& $91.6$
& $126.9$ ($\boldsymbol{-71\%}$) \\

\multirow{2}{*}{$288$}
& \multirow{2}{*}{$5$}
& \multirow{2}{*}{$896.3^{\dagger}$}
& \multirow{2}{*}{$589.8$ ($\boldsymbol{-34\%}$)}
& \multirow{2}{*}{$803.9^{\dagger}$ ($\boldsymbol{-10\%}$)}
& $1\%$
& $107.1$
& $365.5$ ($\boldsymbol{-59\%}$) \\
& & & &
& $5\%$
& $98.6$
& $212.3$ ($\boldsymbol{-76\%}$) \\

\multirow{2}{*}{$144$}
& \multirow{2}{*}{$10$}
& \multirow{2}{*}{$964.6^{\dagger}$}
& \multirow{2}{*}{$663.6^{\dagger}$ ($\boldsymbol{-31\%}$)}
& \multirow{2}{*}{$798.3^{\dagger}$ ($\boldsymbol{-17\%}$)}
& $1\%$
& $95.1$
& $476.5$ ($\boldsymbol{-51\%}$) \\
& & & &
& $5\%$
& $93.8$
& $282.1$ ($\boldsymbol{-71\%}$) \\

\multirow{2}{*}{$288$}
& \multirow{2}{*}{$10$}
& \multirow{2}{*}{$1982.9^{\dagger}$}
& \multirow{2}{*}{$1096.7^{\dagger}$ ($\boldsymbol{-45\%}$)}
& \multirow{2}{*}{$1645.6^{\dagger}$ ($\boldsymbol{-17\%}$)}
& $1\%$
& $119.4$
& $512.4$ ($\boldsymbol{-74\%}$) \\
& & & &
& $5\%$
& $111.7$
& $341.0$ ($\boldsymbol{-83\%}$) \\

\multirow{2}{*}{$144$}
& \multirow{2}{*}{$20$}
& \multirow{2}{*}{$2199.1^{\dagger}$}
& \multirow{2}{*}{$1300.1^{\dagger}$ ($\boldsymbol{-41\%}$)}
& \multirow{2}{*}{$1568.4^{\dagger}$ ($\boldsymbol{-29\%}$)}
& $1\%$
& $101.3$
& $573.2$ ($\boldsymbol{-74\%}$) \\
& & & &
& $5\%$
& $100.2$
& $325.4$ ($\boldsymbol{-85\%}$) \\

\hline
\end{tabular}
\caption{
Comparison of the traditional full-scale centralized stochastic MPC scheme (FS-MPC), TSA-MPC, ADMM-MPC, and the proposed Algorithm~\ref{alg:control_algorithm}
in terms of the average runtime per MPC iteration
as both the number of time steps $|\boldsymbol{K}|$ and the number of scenarios $|\boldsymbol{\Omega}|$ in the dispatch problem increase.
For Algorithm~\ref{alg:control_algorithm}, the average number of clusters required per MPC iteration at convergence and the optimality gap threshold $\epsilon^{\mathrm{thr}}$ are also reported.
Relative runtime differences with respect to FS-MPC are highlighted in bold in parentheses,
while average runtimes exceeding the prescribed 10-minute limit for computing control actions are denoted by $\dagger$.
}
\label{tab:comp_res}
\end{table}

Finally, Table~\ref{tab:comp_res} compares the computational performance of the traditional full-scale centralized stochastic MPC scheme introduced in Subsection~\ref{subsec:SMPC},
the proposed temporally aggregated distributed stochastic MPC scheme implementing Algorithm~\ref{alg:control_algorithm},
and two additional benchmark MPC schemes.
The latter comprise a temporally aggregated MPC scheme (TSA-MPC), which employs TSA to reduce the temporal dimension of the dispatch problem without scenario-wise decomposition,
and a distributed MPC scheme (ADMM-MPC), which employs ADMM to decompose the dispatch problem across scenarios without temporal aggregation.
TSA-MPC adopts the same TSA procedure as Algorithm~\ref{alg:control_algorithm} whereas ADMM-MPC employs the ADMM formulation of Subsection~\ref{subsec:DSMPC_TSA} directly on the full-scale dispatch problem.
As established in the preceding section, both TSA-MPC and ADMM-MPC provide, at each MPC iteration, a lower bound on the optimal objective function value of the original full-scale MPC scheme.
The corresponding upper bounds are computed analogously to that of Algorithm~\ref{alg:control_algorithm}.
All parameters of TSA-MPC and ADMM-MPC are set consistently with their counterparts in Algorithm~\ref{alg:control_algorithm},
with the optimality gap threshold fixed at $\epsilon^{\mathrm{thr}}=1\%$.
These two benchmark control schemes enable the performance of Algorithm~\ref{alg:control_algorithm} to be assessed not only against the traditional full-scale MPC scheme,
but also against MPC schemes that independently employ the two complexity-reduction methods investigated in this study, namely TSA and mathematical decomposition.

The reported results in Table~\ref{tab:comp_res} are obtained by simulating the operation of the considered MPC schemes over the first week of 2017. 
For both TSA-MPC and the proposed Algorithm~\ref{alg:control_algorithm},
the average inflow scenario is employed as the clustering feature within the adopted sliding window clustering technique.
The reference power signal is generated by perturbing the daily profile shown in Fig.~\ref{fig:fig9_power_gen_example} 
using multiplicative random factors independently sampled from the interval $[0,2]$ for each time step of the MPC prediction horizon.

As shown in Table~\ref{tab:comp_res}, the computational benefits of TSA-MPC and ADMM-MPC are limited for the smallest dispatch problem considered,
with runtime variations of $-9\%$ and $+2\%$, respectively, relative to the traditional full-scale MPC scheme.
However, their computational gains become increasingly pronounced as the problem size grows,
highlighting the potential of TSA and mathematical decomposition to substantially reduce the computational burden of large-scale dispatch problems.
Consistent with their respective formulations, the computational benefits of TSA-MPC increase primarily with the temporal dimension of the dispatch problem,
whereas those of ADMM-MPC increase with the number of scenarios.

More importantly, the proposed temporally aggregated distributed stochastic MPC scheme implementing Algorithm~\ref{alg:control_algorithm} consistently achieves the lowest average runtime across all problem configurations considered in the table.
For the smallest dispatch problem, it reduces runtime by $29\%$ and $54\%$ for optimality gap thresholds of $1\%$ and $5\%$, respectively;
for the largest instance, these reductions increase to $74\%$ and $85\%$, respectively.
The results also highlight the explicit trade-off between computational effort and certified solution accuracy enabled by Algorithm~\ref{alg:control_algorithm}:
by appropriately tuning the optimality gap threshold, the decision-maker can balance computational effort against the desired level of solution accuracy.

Overall, these results show that, for the considered case study,
integrating TSA and mathematical decomposition effectively reduces the computational burden of nonconvex stochastic CH-vRES dispatch problems along both their temporal and scenario dimensions,
yielding substantial computational advantages over traditional full-scale MPC schemes and solution methods that independently employ either TSA or mathematical decomposition.

Notably, the proposed combination of TSA and mathematical decomposition not only significantly improves the computational efficiency of the resulting temporally aggregated distributed stochastic MPC scheme relative to all benchmarked MPC schemes,
but, more importantly, restores tractability for large-scale instances of the considered stochastic nonconvex dispatch problem, for which the alternative MPC schemes fail to compute control actions within the prescribed 10-minute time limit.

\section{Conclusion and Future Work}
\label{sec:conclusions}

This paper addresses the joint dispatch problem of cascaded run-of-the-river hydropower plants and vRES,
specifically wind and solar photovoltaic units.
The dispatch problem is first formulated as a two-stage stochastic MIQP problem
to be solved via a traditional full-scale centralized stochastic MPC scheme.
To enhance scalability, the problem is temporally aggregated using TSA
and subsequently decomposed across scenarios via consensus ADMM.
The resulting temporally aggregated distributed stochastic MPC scheme, as implemented in the proposed Algorithm~\ref{alg:control_algorithm}, simultaneously reduces the temporal dimension of the original dispatch problem while enabling scenario-wise parallelization.

Remarkably, our main theoretical result in Subsection~\ref{subsec:SMPC_TSA} demonstrates that the proposed TSA method ensures the construction of a controller that consistently yields a lower bound on the optimal objective function value of its full-scale counterpart.
As detailed in Subsection~\ref{subsec:DSMPC_TSA}, this property is retained under scenario-wise decomposition.
Building upon these theoretical results, Algorithm~\ref{alg:control_algorithm} provides rigorously validated upper and lower bounds on the optimal objective function value of the original dispatch problem.
Consequently, the proposed control algorithm is endowed with a formal performance guarantee in the form of a theoretically certified optimality gap, computed as the relative difference between the derived bounds.
Notably, Algorithm~\ref{alg:control_algorithm} yields a feasible solution for the original dispatch problem at every iteration.

The reported numerical results validate the performance of the proposed controller across diverse operational settings, 
spanning both deterministic and stochastic instances of the dispatch problem, varying input data, and different parameter configurations for Algorithm~\ref{alg:control_algorithm}.
Notably, the proposed integration of TSA and mathematical decomposition outperforms not only the traditional full-scale MPC scheme
but also alternative MPC schemes that employ TSA or mathematical decomposition independently, rather than integrating both complexity-reduction methodologies within a unified control scheme.
Specifically, the proposed control scheme reduces computational runtime by up to $74\%$ relative to traditional full-scale MPC when a solution accuracy of at least $99\%$ is required, and by up to $85\%$ when the accuracy requirement is relaxed to $95\%$.
More importantly, the proposed controller shows superior computational scalability as the size of the nonconvex dispatch problem increases, 
restoring computational tractability for problem instances in which all benchmarked control schemes fail to compute control actions within the prescribed time limit.

We remark that a direct application of traditional TSA methods or consensus ADMM to the original dispatch model \eqref{eq:cen_sto_MPC} would, in general, result in a purely heuristic control scheme.
In particular, the presence of nonconvexities prevents the standard convergence guarantees of consensus ADMM from holding,
while the approximation accuracy of traditional TSA methods cannot be rigorously quantified.

In contrast, the proposed integration of TSA and consensus ADMM within Algorithm~\ref{alg:control_algorithm} preserves a rigorous theoretical foundation while simultaneously enhancing computational efficiency.
Specifically, the proposed control algorithm provides \textbf{two key practical benefits} for the decision-maker. 
First, it delivers substantial \textbf{computational savings}, as shown by the reported numerical results, which consistently indicate a significant reduction in computational effort relative to traditional MPC schemes. 
Second, it is formally shown to be \textbf{trustworthy}, as it is equipped with a rigorous performance guarantee that enables the accuracy of the computed control actions to be systematically assessed through the derived objective function bounds.

Future research will focus on extending the proposed methodology to more comprehensive formulations of the CH-vRES dispatch problem, incorporating additional system dynamics, such as water travel times, as well as network constraints and market objectives.
Moreover, tailored TSA methods will be investigated to enhance the controller's computational performance by accelerating the convergence of the derived objective function bounds.

\section*{Acknowledgements}
Funded by the European Union (ERC, NetZero-Opt, 101116212). Views and opinions expressed are however those of the authors only and do not necessarily reflect those of the European Union or the European Research Council. Neither the European Union nor the granting authority can be held responsible for them.

\printcredits

\bibliographystyle{elsarticle-num}

\bibliography{cas-refs}



\end{document}